\PassOptionsToPackage{draft}{hyperref}  
\documentclass[12pt]{article}
\usepackage{listings}
\usepackage{amsmath}
\usepackage{amssymb}
\usepackage{graphicx}
\usepackage{rotating}
\usepackage{setspace}
\usepackage[table]{xcolor}
\usepackage{natbib}
 \bibpunct[, ]{(}{)}{,}{a}{}{,}%
 \def\newblock{\ }%
\usepackage{multirow}
\usepackage[letterpaper]{geometry}
\usepackage[]{authblk}
\usepackage{verbatim}
\usepackage{url}
\usepackage{multirow}
\usepackage{fullpage}
\usepackage{subcaption}
\usepackage{booktabs}
\usepackage{graphics}
\usepackage[inline]{enumitem}
\usepackage{longtable}
\usepackage{pdflscape}
\usepackage{afterpage}
\usepackage{adjustbox}
\usepackage{hyperref}

\usepackage{algpseudocode,algorithm,algorithmicx}

\algrenewcommand\algorithmicrequire{\textbf{Precondition:}}
\algrenewcommand\algorithmicensure{\textbf{Postcondition:}}
\allowdisplaybreaks[4]

\newtheorem{theorem}{Theorem}
\newtheorem{proposition}{Proposition}
\newtheorem{corollary}{Corollary}

\numberwithin{equation}{section}

\usepackage{rotating}
\usepackage{tabularray}
\usepackage{tikz}
\usepackage{pgfplots}
\usetikzlibrary{intersections}
\usepackage{caption}

\definecolor{mygreen}{rgb}{0.0, 0.42, 0.24}
\definecolor{bostonuniversityred}{rgb}{0.8, 0.0, 0.0}

\pgfplotsset{compat=1.18}
\begin{document}
\date{}
\renewcommand\Affilfont{\small}

\title{Inverse Mixed Integer Optimization: An Interior Point Perspective}
\author[1]{Samir Elhedhli\thanks{corresponding author; elhedhli@uwaterloo.ca}}
\author[1]{G\"{o}ksu Ece Okur\thanks{gokur@uwaterloo.ca}}

\affil[1]{Department of Management Science and Engineering,  University of Waterloo\\200 University Avenue West, Waterloo, Ontario, Canada}

\maketitle
\onehalfspacing

\vspace*{-2cm}
\begin{abstract}

We propose a novel solution framework for inverse mixed-integer optimization based on analytic center concepts from interior point methods. We characterize the optimality gap of a given solution, provide structural results, and propose models that can efficiently solve large problems. First, we exploit the property that mixed-integer solutions are primarily interior points that can be modeled as weighted analytic centers with unique weights. We then demonstrate that the optimality of a given solution can be measured relative to an identifiable optimal solution to the linear programming relaxation. We quantify the absolute optimality gap and pose the inverse mixed-integer optimization problem as a bi-level program where the upper-level objective minimizes the norm to a given reference cost, while the lower-level objective minimizes the absolute optimality gap to an optimal linear programming solution. We provide two models that address the discrepancies between the upper and lower-level problems, establish links with noisy and data-driven optimization, and conduct extensive numerical testing. We find that the proposed framework successfully identifies high-quality solutions in rapid computational times. Compared to the state-of-the-art trust region cutting plane method, it achieves optimal cost vectors for 95\% and 68\% of the instances within optimality gaps of $e$-2 and $e$-5, respectively, without sacrificing the relative proximity to the nominal cost vector. To ensure the optimality of the given solution, the proposed approach is complemented by a classical cutting plane method. It is shown to solve instances that the trust region cutting plane method could not successfully solve as well as being in very close proximity to the nominal cost vector.

\end{abstract}

\maketitle
 \newpage

\section{Introduction}\label{sec:intro}
Inverse Optimization is an active field of research that originated in \cite{Ahuja} and prospered in critical application domains such as intensity-modulated radiation therapy \citep{ChanMahmoodZhu}, network flow models \citep{yang1997inverse} and machine learning applications \citep{bertsimas2015data}.
 The basic inverse optimization problem seeks to find parameters, closest to a reference point, that render observable solutions optimal. 
Given a matrix $A$ $ \in \mathbb{R}^{m\times n}$, a right-hand side vector $b$ $ \in \mathbb{R}^{m}$,  an observed feasible solution $\hat{x}$ $ \in \mathbb{R}^{n}$, and a reference cost vector $\mathring{c}$ $ \in \mathbb{R}^{n}$, an inverse optimization problem seeks to find a cost vector $c$, closest to  $\mathring{c}$  in the $\ell_1$ norm for which $\hat{x}$ is optimal to $\min \{c^\top x: Ax=b, x \ge 0\}$. The inverse linear optimization problem can be set up as the following bi-level optimization problem:
\[
\begin{array}{rl}
\text{\bf[ILO($\hat{x},\mathring{c}$)]}: & \\
 \mbox{ min } & ||c-\mathring{c}||_1 \\
\text{  s.t. } &
 \hat{x} \in \text{argmin }\left\{ c^\top x:Ax=b, x\geq 0\right\} . 
	\end{array} 
\]
Using the well-established primal-dual optimality conditions, the bi-level optimization problem is equivalent to the single-level optimization problem:
\[
\begin{array}{rl}
\text{\bf[ILO($\hat{x},\mathring{c}$)]}: & \\ 
\mbox{ min} & ||c-\mathring{c}||_1 \\
\text{ s.t. }& A^\top y+s=c,s \geq 0 \text{ and }\\
& (\hat{x},y,s) \text{ are  primal-dual optimal.}
	\end{array}  \]
Mathematically, primal-dual optimality can be enforced by strong duality:  $b^\top y=c^\top \hat{x}$ or by complementarity:    $S\hat{x}=0$ where $S$ is the diagonal of $s$. 
We will adopt the latter for reasons that will become clear later. In that  case, \text{\bf[ILO($\hat{x},\mathring{c}$)]} is equivalent to:
\[
\begin{array}{rl}
 \min ||c-\mathring{c}||_1 & \\
\text{s.t.} &  a_{.i}^\top y\leq c_{i}, \;\;   i \in { \bar{\cal I}} \\
& a_{.i}^\top y=c_{i}, \;\;  i \in {\cal I}, 
	\end{array} 
	\] 
 where $a_{i.}$ are the rows of $A$,
 ${\cal I}=\left\{i=1,...,n; \hat{x}_i >0\right\}$ and ${\bar{\cal I}}=\left\{i=1,...,n;\hat{x}_i=0\right\}$. 
The mixed integer counter-part of the linear inverse optimization problem \text{\bf[IMIO($\hat{x},\mathring{c}$)]} is:
\[\begin{array}{rll}
\text{\bf[IMIO($\hat{x},\mathring{c}$)]}: & &\\
\mbox{ min} &  ||c-\mathring{c}||_1& \\
\text{ s.t. } & \hat{x} \in \text{argmin }\{ &  c^\top x:Ax=b, x\geq 0, \; \\
&&\text{some } x \text{ are integer}\}.
	\end{array} 
	\]
In the absence of a strong primal-dual theory for integer programs, one can think of the mixed integer point $\hat{x}$ as being  an interior point that does not satisfy complementarity. In other words, there exists a dual feasible solution $(y,s)$ such that $S\hat{x}$ is equal to  $\varepsilon \geq 0$. 
Finding the optimal $\varepsilon$ could be optimized in a bi-level optimization problem: 
\[
\begin{array}{rll}
\text{\bf[ILO($\hat{x},\mathring{c}$)]}: && \\
\min      & ||c-\mathring{c}||_1    &                    \notag \\
	\text{s.t.  } \ & (\varepsilon,y,s)=\text{argmin }\{e^\top\varepsilon: 
	& \\&A^\top y+s=c,  s \geq 0, S\hat{x}=\varepsilon\ \}. \notag&
		\end{array}	
		\]
This is the main idea of the proposed approach. We will lay out its theoretical foundation and its relationship to interior point methods and weighted analytic centers, provide optimization models that produce an  $\varepsilon$-optimal solution, and   propose efficient approaches  that can be used  to handle large practical problems.
We first set the inverse optimization problem as one that minimizes the gap between the given mixed integer solution and a  Linear Programming (LP) solution. We prove that the proposed model identifies a cost vector,  an LP optimal solution, and quantifies the absolute optimality gap, which guarantees that  the resulting cost vector renders the given mixed integer point near-optimal.  In testing, though, it would be optimal for most instances.   Additionally, we present a linear programming version of the inverse mixed integer optimization problem that  has a huge impact on solution time.   
Second,  we position  the approach  within the literature, discuss its ties with noisy and data-driven inverse optimization  and provide closed-form solutions similar to  \citet{chan2019inverse}.
Third, we set up the  inverse optimization as  a bi-level program where  the upper-level problem minimizes the norm of the cost vector from the nominal one and the lower-level problem minimizes the absolute optimality gap of the given solution. Due to the dominance of the upper-level objective, the lower-level objective which minimizes the  LP gap, deteriorates considerably. To fix this,  we propose two remedies. First,  a constraint is added to limit the LP gap with respect to the cost deviation norm, which we refer to as the \emph{tolerance model}.
Second, we use a weighted objective that accounts for both the norm and the absolute LP gap, which we refer to as the \emph{ bi-objective model}. Both models provide excellent results as measured by the deviation of $\hat{x}$ from optimality and the deviation of $\hat{c}$ from the reference vector $\mathring{c}$. 
Out of the 177 instances tested and for a relative optimality gap of  $e$-2, 
the tolerance and bi-objective models solve  88\% and  79\% of the instances, respectively and combined they solve 95\% of them. 
These percentages drop to  43\%, 59\%, and 68\% for a relative optimality gap of  $e$-5. 
Overall, the average optimality gaps are 1.3\% and 3.3\%.
Although the norm is not  guaranteed to be optimal, the obtained cost vectors are not far from the nominal cost as their relative norm deviations are on average  37.9\% and 34.9\%.
 Furthermore, the models solve very  fast, with  average CPU times of .43 and .25 seconds, respectively.
Finally, to assert the optimality of  $\hat{x}$, we iterate through  the classical cutting plane method by adding cuts based on feasible solutions generated. 
Testing on the smaller instances, reveals the merit of the approach as it solved all the instances that the state-of-the-art method of \cite{BodurChanZhu} could not solve within their time limit.  \\

 The remainder of the paper is organized as follows.
 In Section \ref{sec:literature}, we review the literature.
In Section \ref{sec:LP}, we present the revised simplex and interior point approaches to  linear inverse optimization. 
 In Section \ref{sec:MIP},
we establish the main theoretical foundation and provide models,  structural results, and ties  to noisy inverse optimization.
 In Section \ref{General}, we discuss the general framework that can handle large practical problems and propose two models, the tolerance and the  bi-objective models.
 In Section \ref{sec:exp}, we test  and compare against the state-of-the-art trust region cutting plane method. In Section \ref{sec:con}, we conclude and provide future directions.\\

{\bf Notation:}
$e$ is the vector of ones of length $n$; $e_{p}$ is a unit vector with the $p$th component being 1. ${\bf x}^{-1}$ denotes the component-wise inverse of ${\bf x}$ while ${\bf X}=$ diag$({\bf x})$. $A$ is $m\times n$ matrix with rows $a_{i.}$ $i=1,...,m$ and columns $a_{.j}$ $j=1,...,n$. 

\section{Literature Review}\label{sec:literature}

Inverse Optimization has been successfully applied to several practical problems including network flow, machine learning, bi-level optimization, and healthcare \citep{ChanMahmoodZhu}. Early work  focused on inferring parameters for specific problems such as network flow \citep[e.g.][]{yang1997inverse}, and on estimating arc costs in shortest path problems \citep{burton1992}.

According to a  recent survey by \cite{ChanMahmoodZhu},  inverse optimization problems can be  classified into three groups: classical inverse optimization, data-driven inverse optimization, and alternative learning paradigms. 
Classical inverse optimization focuses on finding problem parameters that render a given observed solution  optimal, whereas data-driven inverse optimization focuses on finding parameters such that a given observed solution is near optimal as measured by its proximity to optimality.
\citep{chan2019inverse,aswani2018inverse}. Alternative learning paradigms cover topics such as inverse online learning and inverse reinforcement learning.
This work focuses on classical and data-driven inverse mixed-integer optimization. As mixed integer solutions are typically interior, we leverage interior point concepts to  provide a solution methodology. For this reason, we focus the literature review on  solution methodologies for inverse mixed integer optimization.


 Although inverse linear programming  and some special cases have been successfully solved,  inverse mixed integer optimization is still challenging to solve. Existing solution methods have shown promising success, but further advancements are needed to tackle practical large problems successfully.
The first solution methodology is due to 
\cite{Wang2009}, who  devise a cutting plane approach  based on the  $\ell_{1}$ and $\ell_{\infty}$ norms. The approach 
iterates between the solution of a master problem and a subproblem.
The master problem is  a linear program that minimizes the norm deviation and proposes a cost vector that renders the  given solution optimal with respect  to a collection of feasible solutions. The subproblem is  a mixed integer program (MIP) that verifies  the optimality of the cost vector and identifies other feasible solutions to generate further cuts.  
This classical cutting plane approach cuts one or several solutions at a time and may exhibit slow convergence. In addition, each iteration requires the solution of an MIP and  early  termination is not possible \citep{Wang2009}.
To overcome some of these challenges, \cite{duan2011heuristic} propose to parallelize the generation of feasible solutions, and consequently cuts, to reduce the number of iterations and the  computational time.

\cite{schaefer2009inverse} propose two approaches. The first seeks to construct the convex hull of the MIP feasible region  in order to apply  inverse LP, whereas the second approach uses super-additive duality to find a strong dual of the  MIP. Both approaches are computationally burdensome, requiring a large number of iterations  and lacking scalability.
%
\cite{lamperski2015polyhedral} further explore  super-additive duals and consider integer and continuous variables separately. They suggest different value functions for the linear and integer parts and consider the maximum value of their sum as the value function of MIP. If  the MIP has  no continuous variables, their model is equivalent to \cite{schaefer2009inverse}, and if there are no integer variables, their approach is equivalent to \cite{Ahuja}. 

The most recent, and state-of-the-art, cutting plane methodology is due to 
\cite{BodurChanZhu}. They  present a  cutting plane framework with  a cut generation subproblem that returns interior forward-feasible points within  a trust region. The search  is confined to a trust region centered at the given solution  $\hat{x}$ and  only includes points whose $\ell_{1}$ distance to $\hat{x}$ does not exceed a predetermined threshold. Similar to \cite{Wang2009}, they iterate between solving an LP master problem and a trust-region  MIP. If violated cuts are identified, they are added to the master problem. If not,  the sub-region, which is the intersection of the feasible region and the trust region,  is increased. To ensure convergence, the trust region is  removed in select iterations. 
To enhance the approach, the  authors  develop an early stop heuristic and impose a time limit on the cut generation subroutine. When the time limit is reached, an  early stop heuristic returns the point with the maximum violation. If not,  the first feasible violated point encountered after the time limit is returned.
They also suggest a dimension reduction step for cut generation by allowing only a random set of indices of the new point to be different from $\hat{x}$. 
For testing, the authors construct a test bed based on 73 problems from the MIPLIB 2017 library \citep{gleixner2021miplib} and generate 3 optimal solutions for each instance. They compare  the classical cutting plane algorithm of \cite{Wang2009}, the cutting plane with early stop, the cutting plane with trust region, the cutting plane with trust region and early stop, the cutting plane with early stop, trust region and dimension reduction  based on the number of instances  solved within an hour. It is found  that the cutting plane with early stop
 and 
the cutting plane with early stop, trust region and dimension reduction 
are able to solve more instances and have faster solution times. The authors also note that 83 instances remain unsolved by any of the suggested methods.

\section{An interior Point Approach to Inverse Optimization} \label{sec:LP}

To lay the groundwork for the  inverse mixed optimization approach, let us first revisit the linear case and offer a different perspective, based on the two leading solution paradigms for linear programming:  revised simplex and interior point methods.

\subsection{A Revised Simplex Approach to Inverse Linear Optimization}
We follow the same notation as in \cite{BodurChanZhu} and consider the following forward linear optimization problem:
\[\begin{array}{rl}
\text{\bf [FLO(c)]}:&  \min\limits_x\; c^\top x \\
&\text{s.t.} \; Ax=b,  \;  x\geq 0. 
	\end{array}	
\]
which is an LP in standard format with feasible set ${\cal X}=\{x\geq 0: Ax=b\}$.
For a given feasible solution $\hat{x}$ and a reference cost vector $\mathring{c}$, the corresponding inverse optimization problem is:
\[
\begin{array}{rl}
\text{\bf[ILO($\hat{x},\mathring{c}$)]}: &
\min\limits_c\;\;   \|c-\mathring{c}\|_1 \\
&\text{s.t.}\quad \hat{x} \text{ is optimal to }\text{\bf [FLO(c)]},
\end{array}
\]
which is, in fact, a bi-level programming problem of the form:
\[
\begin{array}{rl}
\text{\bf[ILO($\hat{x},\mathring{c}$)]}:
\min\limits_c\quad &  \|c-\mathring{c}\|_1  \\ 
\text{s.t.}& \hat{x} \in \text{argmin} \{c^\top x:Ax=b,x\geq 0\}
\end{array}
\]
The straightforward approach to solve the previous problem is to exploit LP duality or complementarity slackness and cast it as a single-level LP:
\[
\begin{array}{rl}
\text{\bf[ILO($\hat{x},\mathring{c}$)]}:\quad\min\limits_c\quad &  \|c-\mathring{c}\|_1 \\ 
\text{s.t.}\quad& A^\top y+s=c \\
&c^\top\hat{x}=b^\top y \\
&s\geq 0.
\end{array}
\]

For this LP to yield meaningful results, $\hat{x}$ should be on the boundary of ${\cal X}$. Two cases arise. If $\hat{x}$ is an extreme non-degenerate point, i.e. a basic feasible solution, one can deduce the optimal basis ${\bf B}$ given by the columns of A corresponding to $\hat{x}>0$, denoted by ${\bf B}=A_{{\hat x}>0}$. The remaining columns of $A$ constitute the non-basic coefficient matrix ${\bf N}=A_{{\hat x}=0}$. When the solution is degenerate, one could find the basis iteratively by adding columns of A that $\hat{x}_{i} = 0$ until an invertible basis is obtained. Accordingly, let us partition the cost vector into $c_{\bf B}$ and $c_{\bf N}$. For $\hat{x}$ to be optimal to $\text{\bf [FLO(c)]}$, the reduced costs should satisfy:
\[
c^\top_{\bf N}-c^\top_{\bf B} {\bf B}^{-1}{\bf N} \geq 0.
\]
As a result, the inverse linear optimization problem reduces to:
\[
\begin{array}{rl}
\text{\bf[ILO($\hat{x},\mathring{c}$)]}: &\min\limits_c\;\; \;  \|c-\mathring{c}\|_1  \\
&\text{s.t.  }  c^\top_{\bf N}-c^\top_{\bf B}{\bf B}^{-1}{\bf N} \geq 0.
\end{array}
\]
This eliminates the need for dual variables $y,$ and is of smaller dimension.
If $\hat{x}$ is not an extreme point, we  resort to a slightly different formulation that is inspired by the rich field of interior point methods. Using complementarity in lieu of strong duality, we write the optimization problem as:
 \[
\begin{array}{rl}
\text{\bf[ILO($\hat{x},\mathring{c}$)]}:\quad\min\limits_c\quad &  \|c-\mathring{c}\|_1 \\ 
\text{s.t.}\quad& A^\top y+s=c \\
&\hat{x}\odot s = 0\\
&s\geq 0.
\end{array}
\] 
where $\odot$ denotes the component-wise product, aka Hadamart product. Taking advantage of the sign of $\hat{x}$, it reduces to: 
\[
\begin{array}{rl}
\min\limits_c & \|c_{\bf \tilde{B}}-\mathring{c}_{\bf \tilde{B}}\|_1 + \|c_{\bf \tilde{N}}-\mathring{c}_{\bf \tilde{N}}\|_1 \\
\text{s.t.}&  {\bf \tilde{B}}^\top y=c_{\bf \tilde{B}},
\quad {\bf \tilde{N}} ^\top y\leq c_{\bf \tilde{N}},
\end{array}
\]
where ${\bf \tilde{B}}=A_{\{{\hat x}>0\}}$ and ${\bf \tilde{N}}=A_{\{{\hat x}=0\}}$.
It could be solved  as is, or further reduced to 
\[\begin{array}{rl}
\min\limits_c  & \|c_{\bf \tilde{B}}-\mathring{c}_{\bf \tilde{B}}\|_1 + \|c_{\bf \tilde{N}}-\mathring{c}_{\bf \tilde{N}}\|_1 \\
\text{s.t.}&
\;\;  {\bf \tilde{N}} ^\top [{\bf \tilde{B}}{\bf \tilde{B}}^\top]^{-1} {\bf \tilde{B}} c_{\bf \tilde{B}} \leq c_{\bf \tilde{N}}  
\end{array}
\]
for $[{\bf \tilde{B}}{\bf \tilde{B}}^\top]$ invertible.

Now let us consider the case where $\hat{x}$ is interior to ${\cal X}$, but expected to be  close to the boundary. The inverse optimization problem then reduces to finding  a cost vector so that $\hat{x}$ is near optimal, i.e., within an $\varepsilon$  of the optimal solution. This is closely related to noisy inverse optimization of \cite{chan2014generalized} and \cite{chan2019inverse}, which is discussed   in the following section.

\subsection{Noisy Inverse Optimization}

\cite{chan2019inverse} focus on estimating cost parameters in noisy inverse optimization when the provided solution may not be exactly optimal. In such cases, applying exact  inverse optimization methods may lead to infeasibility or a zero-cost vector. To address this, the authors minimize the deviation from optimality  by introducing various error measures and provide closed-form solutions under a normalized $\ell_{1}$ norm, i.e. $\|c\|_1$=1.
They introduce three error functions, starting with the Generalized Inverse Optimization model (GIO), which minimizes the distance to a boundary point. This distance can be measured using any norm, and they provide closed-form solutions, reducing the search for an optimal objective coefficient vector to comparing the normalized constraint coefficients. They explore the $\ell_{1}$, $\ell_{2}$ and $\ell_{\infty}$ norms in more detail and offer  closed-form solutions for the optimal error vectors. Next, they discuss the Absolute Duality Gap  model (GIOA), which measures the error in the strong duality constraint. The novelty of this approach is that when the cost coefficients are non-negative, the problem can be formulated as a linear program. The authors establish connections between GIOA and GIO and provide closed-form solutions for GIOA as well. Finally, they consider the Relative Duality Gap  model (GIOR), which minimizes the ratio between primal and dual objectives. Again, they draw connections between GIOR and GIO and  offer closed-form  solutions. They note that when the dual objective is non-negative, this problem can also be solved as a linear program. Additionally, they discuss measures for goodness of fit, noting that structural constraints on the error or coefficient vector may render closed-form solutions inapplicable. For example, structural constraints on $c$ may prevent a direct solution of GIO, making GIOA or GIOR preferable. However, they suggest that the GIO closed-form solution can be adapted when there are constraints on the error.  
Earlier work, 
\cite{chan2014generalized},  focuses on estimating cost parameters in noisy inverse optimization, specifically applied to intensity-modulated radiation therapy. Practical treatments involve balancing conflicting objectives and are typically modeled as a single objective using weights. The work aims  to predict the weight values that make previous treatment plans nearly optimal.

\subsection{An Interior Point Approach to Inverse Linear Optimization}\label{sec:inter-point}

Let us assume that $\hat{x}$ is strictly interior in the sense that $ A\hat{x}=b$ and  $\hat{x}>0$. Such an interior point is the unique minimizer of a  weighted primal potential function, for certain weights  
 \citep{AndersonJibrin,KarimiMoazeniTuncel}:  
\[
\min \;  c^\top x -\sum_{i=1}^{n}\varepsilon_i \ln(x_{i})\;\;
\text{s.t.} \;\; Ax =b, \;x>0.
\]
It has a corresponding dual counterpart $(\hat{y},\hat{s})$, the unique maximizer of the weighted dual potential function: 
 \[
\max\;  b^\top y +\sum_{i=1}^{n}\varepsilon_i \ln(s_i) \;\;
\text{s.t.} \;\;  A^\top y + s = c,\;s>0.
\]
The primal-dual solutions $(\hat{x},\hat{y},\hat{s})$  are interior to their respective feasible regions:
\[{\cal F}_p=\{x: Ax=b,x>0\}, \quad \quad 
{\cal F}_d=\{y: A^\top y < c\},\]
which are assumed bounded  with non-empty interior. 
The solutions $(\hat{x},\hat{y},\hat{s})$ are known as \emph{ weighted analytic centres}, and are the solutions of the following system of equations:
 \[
\begin{array}{rl}
Ax=b, x>0 \\
A^\top y +s =c, s>0 \\
Xs=\varepsilon.
\end{array}
\]
The case where $\varepsilon$ are all the same, i.e., $\varepsilon=\mu e$, leads to the analytic center, and the curve defined by the different values of $\mu>0$  is the central path  \citep{RoosTerlakyVial,Ye}.
The reverse is also true. For,  any point $\hat{x}$ interior to ${\cal X}$,  there exists $(y,s,\varepsilon)$ such that 
$A^\top y +s =c, s>0$ and $\hat{X}s=\varepsilon$, which can be reduced to:  
 \[
\begin{array}{rl}
A^\top y +\hat{X}^{-1}\varepsilon =c,  \varepsilon \geq 0.
\end{array}
\]
Furthermore, the primal-dual solution $(x,y,s)$ has an optimality gap of 
\[\begin{array}{rl}
b^\top y-c^\top \hat{x} &= b^\top y-\hat{x}^\top (A^\top y+s) \\
&=b^\top y-y^\top A\hat{x}+\hat{x}^\top s=e^\top\varepsilon.
\end{array}\]

For any given $c$, the minimum gap is achieved by solving:
\[
\min\;\;  e^{\top}\varepsilon \;\;
\text{s.t.}\;\; A^\top y+ \hat{X}^{-1} \varepsilon=\mathring{c}, \;\;\varepsilon \geq 0,
\]
which is a linear program, whose solution  ($\hat{y},\hat{\varepsilon}$) 
guarantees that $\hat{x}$ is within $e^{T}\hat{\varepsilon}$ of optimality. Specifically for $c=\mathring{c}$,  $A^\top \hat{y}+ \hat{s} = \mathring{c}$  with $\hat{s} = \hat{X}^{-1} \hat{\varepsilon}$, one could construct   different $c$ vectors by noting that 
$A^\top \hat{y}+ (\hat{s} - \sigma) =\mathring{c} - \sigma$
for $0 < \sigma < \hat{s}$. This implies that $\hat{x}$ is $\sum_{i=1}^{n} 
(\hat{\varepsilon_{i}}- \sigma_{i}\hat{x}_{i})$-optimal with respect to $c=\mathring{c}-\sigma$. This clearly demonstrates the trade-off between proximity to $\mathring{c}$ and near-optimality.
We leave the most general case where $\hat{x}$ is not strictly interior to the following section on inverse mixed integer optimization.

\section{Inverse Mixed Integer Optimization}\label{sec:MIP}

In this section, we consider the general case where $\hat{x}$ is a mixed integer solution to the forward mixed integer optimization problem.
\[
\begin{array}{rl}
[\text{\bf FMIO(c)}]:& \\
\min\limits_x\quad  c^\top x &\\ 
\text{s.t.}\quad Ax=b&\\ 
 x_j\geq 0 &\quad j \in  {\cal J} =\{1,...,q\} \\ 
 x_j\geq 0 &\quad\text{integer},  j \in \bar{{\cal J}}=\{q+1,...,n\}.  
\end{array}
\]
Its LP relaxation is $\text{\bf [FLO(c)]}$, which we assume feasible and has a finite solution $x^{LP}.$
For a given mixed integer feasible solution $\hat{x}$ and a reference cost vector $\mathring{c}$, the corresponding inverse mixed integer optimization problem is:
\[
\begin{array}{rl}
[\text{\bf IMIO($\hat{x},\mathring{c}$)}]:\;\min\limits_c\quad &  \|c-\mathring{c}\|_1 \\ 
\text{s.t.}\;& \hat{x} \text{ optimal to }[\text{\bf FMIO(c)}].
\end{array}
\]
It is equivalent to:
\[
\begin{array}{rl}
[\text{\bf ICP}(\hat{x},\mathring{c},{\cal H})]:\;\min\limits_c\quad &  \|c-\mathring{c}\|_1 \\ 
\text{s.t.}\quad&c^\top(\hat{x} - x^{h}) \leq 0, \;  x^{h} \in {\cal H},
\end{array}
\]
where ${\cal H}$ 
is the set of extreme points of the convex hull of 
$\{x\geq 0:Ax=b;  x_j$ integer $ j=q+1,...,n$. 
This is the basis of  the cutting plane algorithm of \cite{Wang2009} where the constraints are generated iteratively by first solving $[\text{\bf ICP}(\hat{x},\mathring{c},\overline{\cal H})]$ for a subset of extreme points denoted by ${ \overline{\cal H}} \subseteq {\cal H}$, obtaining a candidate cost vector $\tilde{c}$, and solving  the subproblem:
\[
[\text{\bf FCP}(\tilde{c})]:\quad \min\limits_{x\in{\cal X}} \quad  \tilde{c}^\top x. 
\]
If $\hat{x}$ is an optimal solution to $[\text{\bf FCP}(\tilde{c})]$, the algorithm stops and 
$\tilde{c}$ is the solution to $[\text{\bf ICP}(\hat{x},\mathring{c},{\cal H})]$. If not, the  
solution is added to ${\cal H}$ and a new iteration is executed.
\cite{BodurChanZhu} restrict $\cal X$ to a sub-region that intersects with a trust region centered at $\hat{x}$, which decreases the cut generation time.
 
For our proposed approach, let us define the inverse mixed integer optimization problem  as:
\[\begin{array}{rl}
[\text{\bf IMIO($\hat{x},\mathring{c}$)}]: & \\
\min\limits_c &\|c-\mathring{c}\|_1 \\
\text{s.t.}&  c^\top \hat{x} \text{ is closest to }c^\top x^{LP}.
\end{array}\]
This problem setup is motivated by the fact that an interior optimal solution is the one with the minimum deviation from an  LP relaxation solution as measured by the absolute objective differences.
For a given $c$, the lower-level problem is 
\[\begin{array}{rl}
\min& \sum_{i \in {\cal I}}\varepsilon_{i} \\
\text{s.t.}\ &  A^\top y+s=c \\
& \hat{x}_i s_i = \varepsilon_{i}; \; i \in {\cal I}, \;\;
s\geq 0, \varepsilon \geq 0.
\end{array}
\]
It is, in turn, equivalent to:
\[\begin{array}{rl}
\min\; &\sum_{i \in {\cal I}}\varepsilon_{i}  \\
\text{s.t.}& a_{.i}^\top y+\frac{\varepsilon_{i}}{\hat{x}_{i}}=c_i, i \in {\cal I} \\
& a_{.i}^\top y+s_i=c_i, i \in  {{\bar{\cal I}}}, \;\;
s\geq 0, \varepsilon \geq 0,
\end{array}
\]
and the corresponding primal-dual optimality conditions are:

\begin{align}
a_{.i}^\top y+\frac{\varepsilon_{i}}{\hat{x}_{i}}=c_i, &\; i \in {\cal I} &\label{op1}\\
a_{.i}^\top y+s_i=c_i, s_i \geq 0, &\; i \in  {{\bar{\cal I}}}& \label{op2}\\
A\delta =0, && \label{op3}\\
\delta_{i} \leq \hat{x}_{i}, & \;  i= 1,...,n &\label{op4}\\
\delta_{i} s_i =0, &\; i \in  {{\bar{\cal I}}} &\label{op5}\\
\varepsilon_{i}( \hat{x}_{i} -\delta_{i})=0,   &\;  i \in {\cal I}. &\label{op6}
\end{align}

We are now ready for the first major result:

\begin{theorem} \label{th:feas}\hspace*{-5pt}{\bf:}
There is at least one cost vector $c$, for which the  system of equations (\ref{op1}-\ref{op6}) possesses a  solution.
\end{theorem}

\noindent {\bf Proof:}
All we have to do is prove that 
\[
\min \; c^\top \delta \;\;
 \text{s.t.} \;\; A\delta=0; \;\;
		 \delta_{i}\leq \hat{x}_i, \; i \in {\cal I};\;
		 \delta_{i}\leq 0, \;i \in {{\bar{\cal I}}}\]

is feasible for a certain cost vector $c$ as the rest follows by complementarity.
Given that $\{x:  Ax=b, x \geq 0\}$ is non-empty and $\hat{x}$ is not an extreme point, there should exist an extreme point $\tilde{x} \neq \hat{x}$ with $A\tilde{x}=b, \tilde{x}>0$. Define $\tilde{\delta}$ as $\hat{x}-\tilde{x}$. It follows that
$A\tilde{\delta}=A\hat{x}-A\tilde{x}=b-b=0$; 
$\hat{x}-\tilde{\delta}=\tilde{x}\geq 0$ which gives $\tilde{\delta}\leq \hat{x}$.
Given the definition of ${\cal I}$, it follows that 
$\tilde{\delta}_i\leq \hat{x}_i, \;\; i \in {\cal I}$ and 
$\tilde{\delta}_i\leq 0, \;\;i \in {{\bar{\cal I}}}$.
%
\hfill $\square$
 \\

We now provide the following important result based on  (\ref{op1}-\ref{op6}).

\begin{theorem}
\label{TH2}\hspace*{-5pt}{\bf:}
Given a solution $(\delta^*,y^*,s^{*},\varepsilon^{*})$ 
to (\ref{op1}-\ref{op6}), an optimal solution to the LP relaxation $\text{\bf [FLO(c)]}$ and its dual  is  given by
$x^{LP}=\hat{x}-\delta^*;y^{LP}=y^*; s_i^{LP}=\frac{\varepsilon_{i}^{*}}{\hat{x}_{i}}, i \in  {\cal I};
s_i^{LP}=s^{*}_{i},  i  \in  \bar{{\cal I}}.$
\end{theorem}

\noindent {\bf Proof:} It suffices to prove that $x^{LP},y^{LP},s^{LP}$ are a primal-dual optimal  to 
$\{\min c^\top x : Ax=b; x\geq 0 \}$ and its dual $\{\max b^\top y : A^\top y+s=c; s\geq 0 \}$ as verified through primal feasibility, dual feasibility and strong duality.
\begin{itemize}
	\item 
Primal feasibility:  $Ax^{LP}=b$ follows from  $A\hat{x}=b$  and $A\delta^*=0$ and $x^{LP} \geq 0$
from $\delta_{i} \leq \hat{x}_{i}$.
\item 
 Dual feasibility: $A^\top y^*+s^*=c$ follows from $a_{.i}^\top y^*+\frac{\varepsilon_{i}^{*}}{\hat{x}_{i}}=c_i, i \in {\cal I} $ and 
$a_{.i}^\top y^*+s^*_i=c_i, s_i \geq 0, i \in  \bar{{\cal I}}$ and obviously $s^{LP}\geq 0$.
\item  Strong duality: $c^\top x^{LP} -b^\top y^{LP}$ follows from \ref{op5}  and  \ref{op6} as 
$
\sum\limits_{i} x_i^{LP} s_i^{LP} = \sum\limits_{i} (\hat{x}^*_i-\delta_{i}^*)s_i^{LP} =0.
$
\hfill $\square$

\end{itemize}

Finally, the following result quantifies the relative optimality gap of $\hat{x}$.
\begin{corollary}\hspace*{-5pt}{\bf:}
$c^\top \hat{x} - c^\top x^{LP} = e^{\top}\varepsilon^*$. 
\end{corollary}

\noindent {\bf Proof:} $c^\top \hat{x} - c^\top x^{LP} = c^\top \hat{x} - b^\top y^{*} = c^\top \hat{x} - \hat{x}^\top A^\top y^{*} =  \hat{x}^\top s^{*} =  \sum_{i \in {\cal I}} \hat{x}_{i}\frac{\varepsilon_{i}^*}{\hat{x}_{i}} = e^{\top}\varepsilon^*$.
\hfill{$\square$}

Theorem \ref{TH2} implies that $x^{LP}$ is the closest LP solution to $\hat{x}$ for a given cost vector $c$ as measured by the LP gap $e^\top \varepsilon$. Still,  $\hat{x}$  may  not necessarily be optimal to $[\text{\bf FMIO(c)}]$. 
In this sense, the proposed framework fits within the {\it noisy and data-driven optimization} paradigm of \citet{chan2019inverse} where model-data fit error is  accounted for in making a given solution approximately optimal. Their general modeling framework  is:
\begin{align}
[\text{\bf GIO}(\hat{x})]: 
      \min \;  &  \|\varepsilon\|_{1} \notag \\
       \text{s.t.}\quad & A^{\top}y + s = c  \notag \\
       &    c^{\top}(\hat{x}-\varepsilon) = b^{\top}y  \notag \\
       &   \|c\|_{1} = 1   \notag \\
       &   s \ge 0.   \notag
    \end{align}
The model-data fit error is measured by the minimum norm of $\varepsilon$ by making $\hat{x}-\varepsilon$ optimal. Other measures were explored including the absolute optimality gap. For  our approach, our modeling framework is:
\begin{align}
[\text{\bf IOP}(\hat{x})]: 
\min \;  & \; e^\top\varepsilon \notag \\
\text{s.t. }  & a_{.i}^\top y+\frac{\varepsilon_{i}}{\hat{x}_{i}}-c_{i}=0  \;\;   i \in {\cal I} \label{IOP1c1}\\
     & a_{.i}^\top y+s_{i}-c_{i}=0  \;\;  i \in  \bar{{\cal I}} \label{IOP1c2}\\  
     &\|c\|_1 = 1    \label{IOP1c4}  \\ 
     & \varepsilon_{i} \geq 0, \: i \in {\cal I}; s_{i} \geq 0, \:\: i \in  \bar{{\cal I.}} \label{IOP1c5}
\end{align}
Inspired by  \citet{chan2019inverse}, we identify  optimal solutions  in closed-form expressions. 
\begin{proposition} \hspace*{-5pt}{\bf:}
An optimal solution to $[IOP(\hat{x})]$ is  given by $y^{*} = \frac{e_{p}}{\|a_{p.}\|_1}$, $c^{*} = \frac{a_{p.}}{\|a_{p.}\|_1} $, $s^{*} = 0$, $\varepsilon^{*} = 0$, where  $a_{p.}$ is any  row of $A$ and $e_{p}$ is a unit vector with the $p$th component being 1.
\end{proposition}

\noindent {\bf Proof:} We first show that the given solution is feasible.
Constraints  (\ref{IOP1c4}) and (\ref{IOP1c5}) are obviously satisfied.
Constraints (\ref{IOP1c1}) and  (\ref{IOP1c2}) follow from plugging in $y^*$ and $c^*$. 
Optimality is easy to verify as $\varepsilon \ge 0$ implies that the best objective achievable  is 0. 
\hfill $\square$

We note that when $\hat{x}$ has  one entry, $\hat{x}_p$, that is zero, 
it can be easily verified that   $y^{*} = 0$, $c^{*} = e_{p} $, $s^{*}_{p} = 1$
$s^{*}_{i} = 0$, $i\neq p$ is an optimal solution.
In the following section, we focus on the case where the constraints of the forward optimization problem are inequalities of the form $Ax\leq b.$ Although they can be turned into standard form $Ax+\sigma=b$ by adding non-negative slacks $\sigma$, it may be necessary to require the cost of the  added  slacks to be zero. 
In that case, the forward optimization problem is 
$ \{\min c^{\top}x\: \: \text{s.t.} \:\: Ax + \sigma = b, x,\sigma  \geq 0\}$
and the inverse problem is:  
\begin{align}
[\text{\bf IOP2}(\hat{x}, \hat{\sigma} )]: 
	  \min  & \:  e^\top\varepsilon + e^\top \varepsilon^\sigma \notag \\
  \text{s.t.} \; &a_{.i}^\top y+\frac{\varepsilon_{i}}{\hat{x}_{i}}-c_{i}=0, \;\;    i \in {\cal I} \label{IOP2_c1} \\
   &   a_{.i}^\top y+s_{i}-c_{i}=0, \;\; i \in {\bar{\cal I}} \label{IOP2_c2} \\  
&    y_{i}+\frac{\varepsilon^\sigma_{i}}{\hat{\sigma}_i}=0,  \;\;  i \in {\cal I_{\sigma}} \label{IOP2_c3} \\
   &  y_{i} +\sigma_{i}=0, \;\; i \in {\bar{\cal I}_{\sigma}} \label{IOP2_c4} \\  
   &  \|c\|_1 = 1    \label{IOP2_c6} \\ 
   &   \varepsilon_{i} \geq 0, \;\;   i \in {\cal I}; \varepsilon^\sigma_{i} \geq 0, \;\; i \in {\cal I}_{\sigma} \label{IOP2_cc7} \\
	   &  
	s_{i} \geq 0, \:;\; i \in {\bar{\cal I}}; \sigma \geq 0, \;\;   i \in {\bar{\cal I}_{\sigma}}, \label{IOP2_c7}
\end{align}
where ${\cal I}_{\sigma}=\{i: \hat{\sigma}_i >0, i=1,...,m;\}$  and  $\bar{\cal I}_{\sigma}=\{i:\hat{\sigma}_i=0,i=1,...,m;\}.$
We characterize the optimal solution to $[\text{\bf IOP2}(\hat{x}, \hat{\sigma} )]$ under two cases depending on whether ($\hat{x}, \hat{\sigma} $) is strictly interior or on the boundary.

\begin{proposition}\hspace*{-5pt}{\bf:}
For  $(\hat{x}, \hat{\sigma}) \not> 0$:
\begin{itemize}
	\item 
An optimal solution to $[\text{\bf IOP2}(\hat{x}, \hat{\sigma} )]$, where $\hat{x}$ contains at least one 0 is  given by: \\ $y^{*} = 0$, $c^{*} = e_{p} $, $s^{*}_{i} = 1 $ if $ i = p$; 0 o.w., $\varepsilon^{*} = 0$, where $p$ is the index of any 0 variable in $\hat{x}$.
\item
An optimal solution to $[\text{\bf IOP2}(\hat{x}, \hat{\sigma} )]$, where $\hat{\sigma}$ contains at least one 0 is  given by:\\ $y^{*} = \frac{-1}{\|a_{p.}\|_1} $ if $ i = t$; 0 o.w., $c^{*} = \frac{a_{p.}}{\|a_{p.}\|_1} $, $s^{*}_{i} = \frac{1}{\|a_{p.}\|_1} $ if $ i = p$; 0 o.w., $\varepsilon^{*} = 0$, where $p$ is the index of any 0 variable in $\hat{\sigma}$.
\end{itemize}
\end{proposition}

For these solutions,  the optimal objective is 0. This is different when 
($\hat{x}, \hat{\sigma} $) is strictly interior.

\begin{proposition}\label{theorem1}
For  $(\hat{x}, \hat{\sigma}) > 0$, $[\text{\bf IOP2}(\hat{x}, \hat{\sigma} )]$
\begin{itemize}
	\item[a.] has an  optimal objective  of $\min\{\min_{i \in {\cal I}}\{\hat{x}_{i}\},\min_{i \in {\cal I}_{\sigma}}\{\frac{\hat{\sigma}_i}{\|a_{i.}\|_1}\}\}$.  
	\item[b.] has an optimal solution of  
$y^{*} = \frac{-e_{p}}{\|a_{p.}\|_1}$, $c^{*} = \frac{-a_{p.}}{\|a_{p.}\|_1} $, 
$s^{*} = 0$, $\varepsilon^{*} = 0, \varepsilon^{\sigma*}_{p} = \frac{\hat{\sigma}_{p}}{\|a_{p.}\|_1},
 \varepsilon^{\sigma *} = 0, \, i \in {\cal I}_{\sigma} \setminus \{p\}$ where $p=argmin_{i \in {\cal I}_{\sigma}}\{\frac{\hat{\sigma}_i}{\|a_{i.}\|_1}\}$,  
whenever $\min_{i \in {\cal I}}\{\hat{x}_{i}\} \ge \min_{i \in {\cal I}_{\sigma}}\{\frac{\hat{\sigma}_i}{\|a_{i.}\|_1}\}$.
\item[c.]  has an optimal solution of 
$y^{*} = 0$, $c^{*} = e_{p} $, $\varepsilon^{*}_{s} = 0, \varepsilon^{*}_{p} = \hat{x}_{p}, \varepsilon^{*}_{i} = 0, \,  i \in {\cal I} \setminus \{p\}$, where $p=argmin_{i \in {\cal I}}\{\hat{x}_{i}\}$
whenever $\min_{i \in {\cal I}}\{\hat{x}_{i}\} \le \min_{i \in {\cal I}_{\sigma}}\{\frac{\hat{\sigma}_i}{\|a_{i.}\|_1}\}$.
\end{itemize}
\end{proposition} 

\noindent {\bf Proof:} We first note that $\bar{\cal I}_{\sigma}$ and $\bar{\cal I}$ are both empty, so constraints (\ref{IOP2_c2}) and (\ref{IOP2_c4}) are non-existent and 
${\cal I}=\{1,...,n\}$, ${\cal I}_{\sigma}=\{1,...,m\}$.
Verifying  the  feasibility of the given solutions with respect to the rest of the constraints 
(\ref{IOP2_c1},\ref{IOP2_c3} ,\ref{IOP2_c6}  and \ref{IOP2_c7})  
is straightforward. 
To show that they are optimal under their respective conditions, let us discuss the two cases separately and proceed by a proof by contradiction.\\
{\bf Case 1: $\min_{i \in {\cal I}}\{\hat{x}_{i}\} \ge \min_{i \in {\cal I}_{\sigma}}\{\frac{\hat{\sigma}_i}{\|a_{i.}\|_1}\}$} \\
The objective corresponding to the  suggested solution is $\frac{\hat{\sigma}_p}{\|a_{p.}\|_1}.$
Let us assume that there is a feasible solution 
$(\tilde{y},\tilde{\varepsilon},\tilde{\varepsilon}^\sigma,\tilde{c})$
that has a better objective. Therefore,
\[
 \sum_{i \in {\cal I}}\tilde{\varepsilon}_{i} + \sum_{i \in {\cal I}_{\sigma}}\tilde{\varepsilon}^\sigma_{i} \
<\frac{\hat{\sigma}_p}{\|a_{p.}\|_1}.\]
Constraints (\ref{IOP2_c1}) and (\ref{IOP2_c3}) and the fact that $\hat{x},\hat{\sigma}>0$ imply that
 \[
 \sum_{i \in {\cal I}} \hat{x}_i (\tilde{c}_{i}-a_{.i}^\top \tilde{y}) 
-  \sum_{i \in {\cal I}_{\sigma}} \hat{\sigma}_i \tilde{y}_i 
< \frac{\hat{\sigma}_p}{\|a_{p.}\|_1}
\]
Using 
\[
\hat{x}_i \geq \min_{i \in {\cal I}}\hat{x}_{i} \ge \min_{i \in {\cal I}_{\sigma}}\{\frac{\hat{\sigma}_i}{\|a_{i.}\|_1}\}=\frac{\hat{\sigma}_p}{\|a_{p.}\|_1},
\]
we get
\[
 \frac{\hat{\sigma}_p}{\|a_{p.}\|_1}
 \sum_{i \in {\cal I}}  (\tilde{c}_{i}-a_{.i}^\top \tilde{y}) 
-  \sum_{i \in {\cal I}_{\sigma}} \hat{\sigma}_i \tilde{y}_i 
< \frac{\hat{\sigma}_p}{\|a_{p.}\|_1},
\]
which implies that (note that  $\|a_{p.}\|_1 \neq 0$ and $\hat{\sigma}_p>0$): 
\[
 \sum_{i \in {\cal I}}  (\tilde{c}_{i}-a_{.i}^\top \tilde{y}) 
-  \sum_{i \in {\cal I}_{\sigma}} \frac{\|a_{p.}\|_1}{\hat{\sigma}_p}\hat{\sigma}_i \tilde{y}_i 
< 1.
\]
As $\tilde{\varepsilon},\tilde{\varepsilon}^\sigma, \hat{x}, \hat{\sigma}>0$, it follows that $\tilde{c}_{i}-a_{.i}^\top \tilde{y} >0$ and $|\tilde{c}_{i}-a_{.i}^\top \tilde{y}|=\tilde{c}_{i}-a_{.i}^\top \tilde{y}, i \in {\cal I}$. Similarly, $|\tilde{y}_i|=-\tilde{y}_i, i \in {\cal I_{\sigma}}.$ Therefore, we get:
\[
 \sum_{i \in {\cal I}}  |\tilde{c}_{i}-a_{.i}^\top \tilde{y}| 
+  \sum_{i \in {\cal I}_{\sigma}} \frac{\|a_{p.}\|_1}{\hat{\sigma}_p}\hat{\sigma}_i |\tilde{y}_i| 
< 1.
\]
Using the triangular inequality, $|\tilde{c}_{i}| \leq |\tilde{c}_{i}-a_{.i}^\top \tilde{y}| + |a_{.i}^\top y|$,
and 
$\frac{\|a_{p.}\|_1}{\hat{\sigma}_p} \geq \frac{\|a_{i.}\|_1}{\hat{\sigma}_i}$,
 it follows that:
\[
 \sum_{i \in {\cal I}}  |\tilde{c}_{i}|  - \sum_{i \in {\cal I}}  |a_{.i}^\top \tilde{y}| 
+  \sum_{i \in {\cal I}_{\sigma}} \|a_{i.}\|_1 |\tilde{y}_i| 
< 1.
\]
Finally, given that 
\[\begin{array}{l}
\sum_{i \in {\cal I}_{\sigma}} \|a_{i.}\|_1 |\tilde{y}_i|
=\sum_{i=1}^m \sum_{j=1}^n |a_{ij}\tilde{y}_i| \\
=\sum_{j=1}^n\sum_{i=1}^m  |a_{ij}\tilde{y}_i| 
\geq \sum_{j=1}^n |\sum_{i=1}^m  a_{ij}\tilde{y}_i|\\
\geq\sum_{j=1}^n|a_{.j}^\top \tilde{y}| 
= \sum_{i \in {\cal I}}|a_{.i}^\top \tilde{y}|.
\end{array}\]
 It follows that $\| \tilde{c}\|_1 <1$, which contradicts  (\ref{IOP2_c6}). \\
{\bf Case 2: $\min_{i \in {\cal I}}\{\hat{x}_{i}\} \le \min_{i \in {\cal I}_{\sigma}}\{\frac{\hat{\sigma}_i}{\|a_{i.}\|_1}\}$} \\
The objective corresponding to the  suggested solution is $\hat{x}_p.$
Let us assume that there is a feasible solution 
$(\tilde{y},\tilde{\varepsilon},\tilde{\varepsilon}^\sigma,\tilde{c})$
that has a better objective. Therefore,
\[
 \sum_{i \in {\cal I}}\tilde{\varepsilon}_{i} + \sum_{i \in {\cal I}_{\sigma}}\tilde{\varepsilon}^\sigma_{i} \
<\hat{x}_p\]
Constraints (\ref{IOP2_c1}) and (\ref{IOP2_c3}) and the fact that $\hat{x},\hat{\sigma}>0$ imply that
 \[
 \sum_{i \in {\cal I}} \hat{x}_i (\tilde{c}_{i}-a_{.i}^\top \tilde{y}) 
-  \sum_{i \in {\cal I}_{\sigma}} \hat{\sigma}_i \tilde{y}_i 
< \hat{x}_p
\]
Using 
\[
\hat{x}_i \geq \min_{i \in {\cal I}}\hat{x}_{i} = \hat{x}_p, \]

we get
\[
\hat{x}_p
 \sum_{i \in {\cal I}}  (\tilde{c}_{i}-a_{.i}^\top \tilde{y}) 
-  \sum_{i \in {\cal I}_{\sigma}} \hat{\sigma}_i \tilde{y}_i 
< \hat{x}_p,
\]
which implies that (note that $\hat{x}_p>0$): 
\[
 \sum_{i \in {\cal I}}  (\tilde{c}_{i}-a_{.i}^\top \tilde{y}) 
-  \sum_{i \in {\cal I}_{\sigma}} \frac{{\sigma}_i}{\hat{x}_p} \tilde{y}_i 
< 1.
\]
As $\tilde{\varepsilon},\tilde{\varepsilon}^\sigma, \hat{x}, \hat{\sigma}>0$, it follows that $\tilde{c}_{i}-a_{.i}^\top \tilde{y} >0$ and $|\tilde{c}_{i}-a_{.i}^\top \tilde{y}|=\tilde{c}_{i}-a_{.i}^\top \tilde{y}, i \in {\cal I}$. Similarly, $|\tilde{y}_i|=-\tilde{y}_i, i \in {\cal I_{\sigma}}.$ Therefore, we get:
\[
 \sum_{i \in {\cal I}}  |\tilde{c}_{i}-a_{.i}^\top \tilde{y}| 
+  \sum_{i \in {\cal I}_{\sigma}} \frac{{\sigma}_i}{\hat{x}_p} |\tilde{y}_i| 
< 1.
\]
Using the triangular inequality, $|\tilde{c}_{i}| \leq |\tilde{c}_{i}-a_{.i}^\top \tilde{y}| + |a_{.i}^\top y|$,
and 
$ \hat{x}_p \le \min_{i \in {\cal I}_{\sigma}}\{\frac{\hat{\sigma}_i}{\|a_{i.}\|_1}\} \le \frac{\hat{\sigma}_i}{\|a_{i.}\|_1}$  
%
 it follows that:
\[
 \sum_{i \in {\cal I}}  |\tilde{c}_{i}|  - \sum_{i \in {\cal I}}  |a_{.i}^\top \tilde{y}| 
+  \sum_{i \in {\cal I}_{\sigma}} \|a_{i.}\|_1 |\tilde{y}_i| 
< 1.\]
The rest follows as in the previous case. 
\hfill $\square$

Next, we illustrate on an example. \\
\noindent \textbf{Example:}  Consider the following forward optimization problem
\begin{align}
  \min \: & \:  c^{\top}x \quad \notag \\
   \text{s.t.} &\: -4x_{1}  -3x_{2} \le -19    \notag \\
   &   -x_{1} -3 x_{2} \le -8\: \notag \\
			&
      6x_{1} + x_{2} \le 30  \notag \\
    & -3x_{1} + 5x_{2} \le 17   \notag \\
    &  x_{1}, x_{2}\geq 0  \text{ integer}, \notag
\end{align}
whose feasible region is displayed in the Figure below:
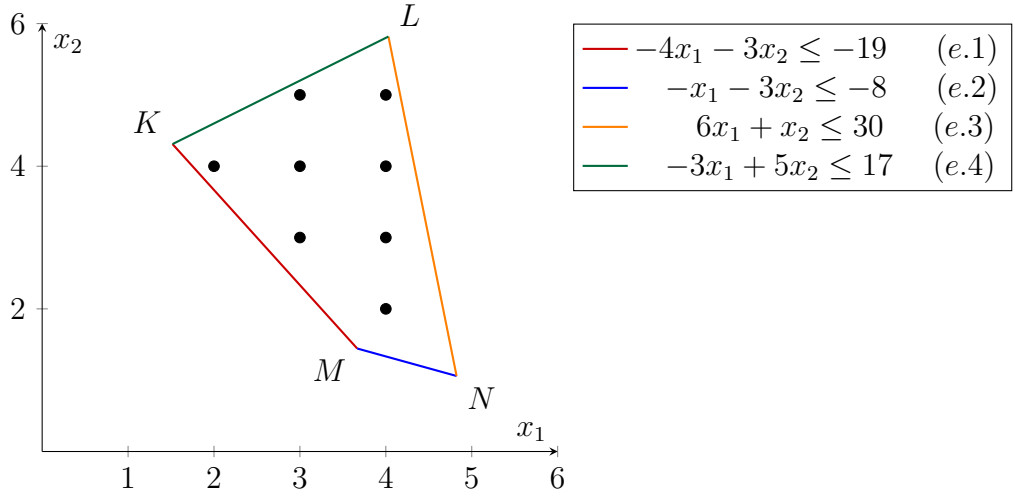
\begin{center}
    \begin{tikzpicture}
        \begin{axis}[
            axis lines = middle,
            xlabel = {$x_1$},
            ylabel = {$x_2$},
            xmin = 0, xmax = 6,
            ymin = 0, ymax = 6,
            legend pos= outer north east,
            clip=false,
            no markers,
            legend cell align=left,
        ]
        
        \addplot[thick, color = bostonuniversityred] coordinates {(44/29,125/29) (11/3,13/9)}; 
        \addlegendentry{$ -4x_1 - 3x_2 \leq -19 \;\;\quad(e.1)$};   
        \addplot[thick, color = blue] coordinates {(11/3,13/9) (82/17,18/17)}; 
        \addlegendentry{$\quad -x_1 - 3x_2 \leq -8 \;\;\quad(e.2)$};     
        \addplot[thick,color = orange] coordinates {(82/17,18/17) (133/33,64/11)}; 
        \addlegendentry{$\quad \quad 6x_1 + x_2 \leq 30 \;\;\quad(e.3)$};       
        \addplot[thick, color = mygreen] coordinates {(133/33,64/11) (44/29,125/29)};
        \addlegendentry{$\quad -3x_1 + 5x_2 \leq 17 \;\quad(e.4)$};

        \addplot[only marks, mark=* ,mark options={fill=black},] coordinates {(2,4)};  
        \addplot[only marks, mark=square*,mark options={fill=black},] coordinates {(3,3)}; 
        \addplot[only marks, mark=triangle*,mark options={fill=black},] coordinates {(3,4)}; 
        \addplot[only marks, mark=* ,mark options={fill=black},] coordinates {(3,5)};  
        \addplot[only marks, mark=square*,mark options={fill=black},] coordinates {(4,2)}; 
        \addplot[only marks, mark=triangle*,mark options={fill=black},] coordinates {(4,3)}; 
        \addplot[only marks, mark=square*,mark options={fill=black},] coordinates {(4,4)}; 
        \addplot[only marks, mark=triangle*,mark options={fill=black},] coordinates {(4,5)}; 

        \node at (axis cs:133/33,64/11) [anchor=south west] {$L$};
        \node at (axis cs:82/17,18/17) [anchor=north west] {$N$};
        \node at (axis cs:11/3,13/9) [anchor=north east] {$M$};
        \node at (axis cs:44/29,125/29) [anchor=south east] {$K$};

        \end{axis}
    \end{tikzpicture}
\end{center}

The solution of  [{\bf IOP2}$(\hat{x}, \hat{\sigma} )$] for  4 different choices of $\hat{x}$ is:
\begin{itemize}
    \item $\hat{x} = (4,2)$ and  $\hat{\sigma} = (3,2,4,19)$, Then using Proposition \ref{theorem1}, the optimal objective of [{\bf IOP2}$(\hat{x}, \hat{\sigma} )$] is min\{min\{4,2\},min\{$\frac{3}{7},\frac{2}{4},\frac{4}{7},\frac{19}{8}\}\}$ = $\frac{3}{7}$, and the optimal $c$ is (0.57, 0.43, 0, 0, 0, 0), which is  the negative of the normalized constraint (e.1) coefficients. For this choice of $c$, however, the  optimal MIP solution is  (2,4).
    \item $\hat{x} = (2,4)$ and  $\hat{\sigma} = (1,4,14,3)$. The optimal objective of [{\bf IOP2}$(\hat{x}, \hat{\sigma} )$] is \\min\{min\{2,4\},min\{$\frac{1}{7},\frac{4}{4},\frac{14}{7},\frac{3}{8}\}\}$ = $\frac{1}{7}$, and the optimal $c$ is (0.57, 0.43, 0, 0, 0, 0), which is the negative of the normalized constraint (e.1) coefficients. For this choice of $c$,  the  optimal MIP solution is $\hat{x}$.
    \item $\hat{x} = (4,5)$ and  $\hat{\sigma} = (12,11,1,4)$. The optimal objective of [{\bf IOP2}$(\hat{x}, \hat{\sigma} )$] is\\ min\{min\{4,5\},min\{$\frac{12}{7},\frac{11}{4},\frac{1}{7},\frac{4}{8}\}\}$ = $\frac{1}{7}$, and the optimal $c$ is (-0.86, -0.14, 0, 0, 0, 0), which is the negative of the normalized constraint (e.3) coefficients. For this choice of $c$,  the  optimal MIP solution is $\hat{x}$.
    \item $\hat{x} = (3,5)$ and  $\hat{\sigma} = (8,10,7,11)$. The optimal objective of [{\bf IOP2}$(\hat{x}, \hat{\sigma} )$] is \\min\{min\{3,5\},min\{$\frac{8}{7},\frac{10}{4},\frac{7}{7},\frac{11}{8}\}\}$ = $\frac{7}{7}$, and the optimal $c$ is (-0.86, -0.14, 0, 0, 0, 0), which is the negative of the normalized constraint (e.3)  coefficients. For this choice of $c$, however, the  optimal MIP solution is (4,5).
\end{itemize}

We note that [{\bf IOP2}$(\hat{x}, \hat{\sigma} )$] does not handle the additional requirement of finding a cost vector that is closest to a reference  cost $\mathring{c}$. In the following section, we address the general case.

\section{The General Approach} \label{General}
Recall that the goal in inverse optimization is to find a cost vector $c$, closest to a reference $\mathring{c}$ that renders $\hat{x}$  optimal or  $\varepsilon$-optimal. In the rest of the work, we focus on the latter and pose the mixed integer inverse optimization problem:
\[
\begin{array}{rl}
[\text{\bf IMIO($\hat{x},\mathring{c}$)}]:\;\min\limits_c\quad &  \|c-\mathring{c}\|_1 \\ 
\text{s.t.}\quad& \hat{x} \:\: \varepsilon\text{-optimal to {\bf [FMIO($c$)]}}.
\end{array}
\]
With the help of (\ref{op1}-\ref{op6}),  it is equivalent to:
\[\begin{array}{rl}
[\text{\bf IMIO($\hat{x},\mathring{c}$)}]: & \notag \\
  \min  &  \sum_{i=1}^{n} (f_{i}+g_{i})  \notag \\
    \text{s.t.}  & a_{.i}^\top y+\frac{\varepsilon_{i}}{\hat{x}_{i}}-f_{i}+g_{i}=\mathring{c}_{i}, \; i \in {\cal I} \notag \\
   & a_{.i}^\top y+s_{i}-f_{i}+g_{i}=\mathring{c}_{i}, \; i \in  \bar{{\cal I}}  \notag \\
   &A\delta = 0 \notag \\
   &\delta \leq \hat{x}  \notag \\
   & \delta_{i} + M z_{i} \geq \hat{x}_{i},  \;   i  \in {\cal I}  \notag \\
   &\varepsilon_{i} + M z_{i} \leq M,\; i \in {\cal I}  \notag \\
   &s_{i} + M z_{i} \leq M,  \;  i \in  \bar{{\cal I}}  \notag \\
   & \delta_{i}\geq -Mz_{i},  \; i \in  \bar{{\cal I}}\notag \\
   & \varepsilon_{i} \geq 0, \: i \in {\cal I};  s_{i} \geq 0, \:\: i \in  \bar{{\cal I}};\notag \\
	& z_{i}\in \{0,1\},f_{i}, g_{i} \geq 0 , \:\:  i= 1,...,n.\notag 
\end{array}
\]
Note that  $(\ref{op5}-\ref{op6})$ are linearized using binary variables $z$ and big $M$
and the  $\ell_{1}$ distance from $c$ to $\mathring{c}$ is linearized using additional variables $f, g$.
Computationally, [\text{\bf IMIO($\hat{x},\mathring{c}$)}] is a binary integer program that may be challenging to solve. The binary variables are introduced to linearize the complementarity slackness conditions. If we use strong duality, instead, and equate objectives of the primal and the dual, we end up with a linear program:
\[
\begin{array}{rl}
[\text{\bf IMIO($\hat{x},\mathring{c}$)}]: & \notag \\
  \min \quad & e^\top f+e^\top g  \notag \\
    \quad \text{s.t.} \quad&  a_{.i}^\top y+\frac{\varepsilon_{i}}{\hat{x}_{i}}-f_{i}+g_{i}=\mathring{c}_{i}, i \in {\cal I} \notag \\
   &a_{.i}^\top y+s_{i}-f_{i}+g_{i}=\mathring{c}_{i}, i \in  \bar{{\cal I}}  \notag \\
   & A\delta = 0  \notag \\
   & \delta \leq \hat{x}  \notag \\
		 &b^\top y+e^\top\varepsilon - \hat{x}^\top f + \hat{x}^\top g=\mathring{c}^\top\hat{x} \notag \\
   & \varepsilon_{i} \geq 0, \: i \in {\cal I};  s_{i} \geq 0, \:\: i \in  \bar{{\cal I}}; \notag \\
	& f_{i}, g_{i} \geq 0 , \:\:  i= 1,...,n. \notag 
\end{array} \]

Note that
$b^\top y+e^\top\varepsilon - \hat{x}^\top f + \hat{x}^\top g=\mathring{c}^\top\hat{x}$, i.e., 
$b^\top y+e^\top\varepsilon=c^\top\hat{x}$  follows from Theorem 2 and  $e^\top\varepsilon=c^\top\delta=c^\top(\hat{x}-x^{LP})=c^\top\hat{x}-c^\top x^{LP}=c^\top\hat{x}-b^\top y.$  Furthermore  $A\delta = 0$  and  $\delta \leq \hat{x}$  could be dropped as the problem decouples, leading to the concise LP:
%
\begin{align}
[\text{\bf IMIO($\hat{x},\mathring{c}$)}]: & \notag \\
\min \quad & e^\top f+e^\top g \notag \\
 \quad \text{s.t.} \quad & a_{.i}^\top y+\frac{\varepsilon_{i}}{\hat{x}_{i}}-f_{i}+g_{i}=\mathring{c}_{i},i \in {\cal I} \label{c11} \\
& a_{.i}^\top y+s_{i}-f_{i}+g_{i}=\mathring{c}_{i},i \in  \bar{{\cal I}}  \label{c12}   \\  
&b^\top y+e^\top\varepsilon - \hat{x}^\top f + \hat{x}^\top g=\mathring{c}^\top\hat{x} \label{c13} \\
&\varepsilon_{i} \geq 0, \: i \in {\cal I};  s_{i} \geq 0, \:\: i \in  \bar{{\cal I}}; \label{cc19} \\
&f_{i}, g_{i} \geq 0 , \:\:  i= 1,...,n. \label{c19}
\end{align}
%
Unfortunately, the bi-level model suffers from  the dominance of the upper-level objective, which, if not carefully addressed, will force the  obvious solution $c=\mathring{c}$. 
 An obvious remedy is to include the lower-level objective alongside the upper-level objective through carefully chosen weights or to link the two objectives. In the following sections, we elaborate on these  two remedies.

\subsection{Weighted Objectives}
First, we  minimize the LP gap alongside the norm in the upper-level objective. For that, we introduce a weight $w$ and  add $w^\top\varepsilon$ to $\|c-\mathring{c}\|_1$, leading to:
\begin{align}
[\text{\bf IMIO($\hat{x},\mathring{c},w$)}]:
	&\:
  \min  \: \sum_{i=1}^{n} (f_{i}+g_{i}) +\sum_{i\in I} w_i \varepsilon_i  \label{obj3}\\
  & \quad \quad \text{s.t.} \: (\ref{c11}) - (\ref{c19}). \notag 
\end{align}
We refer to this model as the \emph{bi-objective model}.
An  interesting choice for $w_i$ is  $\frac{1}{\hat{x}_{i}}$, which mimics constraint $\ref{c11}$.
 We note the following:
\begin{itemize}
	\item For a given solution with $\max(g_i,\varepsilon_i)>0$, other solutions can be constructed by shifting  a fraction of $g_i$ to $\varepsilon_i$ or vice versa.
	\item As $s_i$ is not penalized in the objective, it could safely absorb any values $g_i, i\in \bar{I}$ would assume, eliminating the need for variables  $g_i, i\in \bar{I}$.
	\item When inequality constraints are turned into equality by adding  slack variables, it may be  desirable to set the corresponding cost to 0. This could be achieved  by not including their norm in  $||c-\mathring{c}||_1$, or by  forcing their corresponding $f$ and $g$ to 0. Alternatively, this can be achieved by   making the following transformation once a $c$ vector is identified. Suppose $Ax\leq b$ is turned into $Ax+\sigma=b$, $\sigma\geq0$. Once a cost vector $\hat{c}$ is found, then it follows that 
	$c^\top x=x^\top c^x +\sigma^\top c^\sigma$  
	= $x^\top c^x+(b-Ax)^\top c^\sigma$
	= $[c-A^\top c^\sigma]^\top x+b^\top c^\sigma.$
	Therefore, the new $c$ is set to $c-A^\top c^\sigma$. It may not minimize the norm $||c-\mathring{c}||_1$, though.
\end{itemize}

Let us clarify these notions through the example presented earlier.
Consider $\mathring{c}=(3,1)$ and $\hat{x}= (4,2)$.
Solving \text{\bf IMIO($\hat{x},\mathring{c},w$)} with $w=1$, i.e. 
$( \min  \: \sum_{i=1}^{n} (f_{i}+g_{i}) +\sum_{i\in I} \varepsilon_i )$,
gives the cost vector $\hat{c}=(4/3,1)$ with norm 5/3,
$e^\top \varepsilon^{*}=1$,
an  LP solution  at  point $M$ 
with objective 6.33 (19/3) an  optimal integer solution  $x^{IP}$ at  (2,4) with objective 6.67, $\neq \hat{x}$. 
The full solution is \\
\begin{center}
{\small \singlespacing
\begin{tabular}{|c|c|c|c|c|c|c|c|}\hline
   $\hat{x}$&$x^{IP}$&$x^{LP}$&$\varepsilon^*$&$f^*$&$g^*$&$\hat{c}$& $\mathring{c}$ \\ \hline
   4&   2&   3.67&   0&   0&   1.67&   1.33&   3 \\
   2&   4&   1.44&   0&   0&   0&   1&   1 \\ 
   3&   1&   0&   1&   0&   0&   0&   0 \\
   2&   6&   0&   0&   0&   0&   0&   0 \\
   4&  14&   6.56&   0&   0&   0&   0&   0 \\ 
  19&   3&  20.78&   0&   0&   0&   0&   0 \\ \hline
\end{tabular}	
}\end{center} 

\vspace*{.5cm}
Given that $\varepsilon^{*}_3=1$, we can use constraint (3):  
$y_1+\varepsilon_3/\hat{x}_3-f_3+g_3=0$ with $\hat{x}_3=3$
to shift a portion of $\frac{\varepsilon_3}{\hat{x}_3}$ to $g_3$.
Recall that constraint 3 is
$y_1+\frac{1}{3}-f_3+g_3=0 $.
Shifting $\frac{2}{3}$ of $\frac{\varepsilon_3}{\hat{x}_3}$, which is equal to $\frac{2}{9}$, to $g_3$  gives
$
y_1+\frac{2}{9}-f_3+g_3+\frac{1}{9}=0.
$
The resulting solution is\\
\begin{center}
{\small \singlespacing
\begin{tabular}{|c|c|c|c|c|c|c|c|}\hline
   $\hat{x}$&$x^{IP}$&$x^{LP}$&$\varepsilon^*$&$f^*$&$g^*$&$\hat{c}$& $\mathring{c}$ \\ \hline
   4&   2&   3.67&   0&   0&   1.67&   0.44&   3 \\ 
   2&   4&   1.44&   0&   0&   0&   0.33&   1 \\ 
   3&   1&   0&   0.33&   0&   0.22&   0&   0 \\ 
   2&   6&   0&   0&   0&   0&   0&   0 \\ 
   4&  14&   6.56&   0&   0&   0&   0&   0 \\ 
  19&   3&  20.78&   0&   0&   0&   0&   0 \\\hline
\end{tabular}	
}\end{center} 

\vspace*{.5cm}
This still gives $[2,4]$ as the  optimal solution but
 $\hat{x}=[4,2]$ is now way closer as it  has objective  2.44 as opposed to the optimal 2.22. This comes at the expense of a worse $||\hat{c}-\mathring{c}||_1$ norm of 3.22.

\subsection{Bounding the Optimality Gap}

As mentioned earlier, there is a trade off between the upper-level and  the lower-level objectives and 
solving [\text{\bf IMIO($\hat{x},\mathring{c}$)}] will force $\hat{c}$ to equal
$\mathring{c}$,  causing  $e^{\top}\varepsilon$ to be large and
jeopardizing the optimality of  $\hat{x}$.
One way to force  a small $e^{\top}\varepsilon$ is to put  an upper bound on it. We link such bound to a fraction $\tau$ of the norm  
$||\hat{c}-\mathring{c}||_1$, leading to:   
\begin{align}
[\text{\bf IMIO($\hat{x},\mathring{c},\tau$)}]: \quad
  \min  &  \sum_{i=1}^{n} (f_{i}+g_{i}) \notag \\
   \text{s.t.}  & (\ref{c11})-(\ref{c19})  \notag \\
   & \sum_{i \in {\cal I}} \varepsilon_{i} \le \tau \sum_{i=1,...,n}(f_{i}+g_{i}). \label{ctol}
\end{align}
We refer to this model as the \emph{tolerance model.}
Applied to the example with the same $\mathring{c}$ and $\hat{x}$ and $\tau =1e-3$ leads to the solution:
\vspace*{.5cm}
\begin{center}
{\small \singlespacing
\begin{tabular}{|c|c|c|c|c|c|c|c|}\hline
$\hat{x}$&$x^{IP}$&$x^{LP}$&$\varepsilon^*$&$f^*$&$g^*$&$\hat{c}$& $\mathring{c}$ \\ \hline
  4&2&1.52&0&0&2.995&0.005&3 \\ 
  2&4&4.3&0&0&0.996&0.004&1 \\ 
  3&1&0   &0.004&0&0&0&0 \\ 
  2&6&6.45&0&0&0&0&0 \\ 
  4&14&16.59&0&0&0&0&0 \\ 
 19&3&0&0&0&0&0&0 \\\hline
	\end{tabular}	
}\end{center}

\vspace*{.5cm}
	with $\hat{c}^\top x^{IP}=0.018$, $\hat{c}^\top \hat{x}=0.0195$ and  $\hat{c}^\top x^{LP}=0.0168$. The norm $||\hat{c}-\mathring{c}||_1$ is  3.99, which could be improved by scaling  $\hat{c}$. This gives 
	$\hat{c}=(3;2.25)$ and  $||\hat{c}-\mathring{c}||_1$=1.25,
	$\hat{c}^\top x^{IP}=15$, $\hat{c}^\top\hat{x}=16.5$ and  $\hat{c}^\top x^{LP}=14.26$. If we add the cut, $c^\top\hat{x} \leq c^\top x^{IP}$, increase $\tau$ to 1 and resolve, we get the optimal solution:
	\vspace*{.5cm}
\begin{center}
{\small \singlespacing
\begin{tabular}{|c|c|c|c|c|c|c|c|}\hline
$\hat{x}$&$x^{IP}$&$x^{LP}$&$\varepsilon^*$&$f^*$&$g^*$&$\hat{c}$& $\mathring{c}$ \\ \hline
4&4&3.67&0&0&2&1&3 \\ 
2&2&1.44&0&0&0&1&1 \\ 
3&3&0   &0.67&0&0&0&0 \\ 
2&2&0   &0.22&0&0&0&0 \\ 
4&4&6.56&0&0&0&0&0 \\ 
19&19&20.78&0&0&0&0&0 \\ \hline
	\end{tabular}	
}\end{center}

\vspace*{.5cm}
which makes $\hat{x}$ optimal and gives the  optimal norm of 2. 

\section{Numerical Testing } \label{sec:exp}

To test the viability and efficiency of the proposed approach, we test on  the same MIPLIB 2017 instances as \cite{BodurChanZhu} except for the instances with a lower bound different from zero, which necessitates doubling the number of variables to fit $\{ Ax=b; x\geq 0\}$.  The upper bounds, however,  were incorporated into the coefficient matrix $A$  and the inequality constraints were put in standard format with the help of slack variables. A total of 59 MIP instances were included in our test bed, each with three  $\hat{x}$ solutions,  resulting in  177 instances in total, and denoted by extension $t1$, $t2$, $t3$ as done  in  \cite{BodurChanZhu}.
We divide the instances into ten groups depending on  the number of integer and binary variables with groups 1 to 10 containing 	30 to 100, 130 to 170, 171 to 250, 252 to 379, 396 to 1150,
 1170 to 1372, 1384 to 1914, 2183 to 4150, 4456 to 6642  and 7195 to 11024 binary and integer variables, respectively. 

The input parameters $\hat{x}$ and $\mathring{c}$ were supplied by the authors in \cite{BodurChanZhu}, for a fair comparison against their trust region approach. 
The testing was carried out  on a 
computer with  x64-based Intel(R) Core(TM) i7-4790 processor, 3.60 GHz CPU and 12.0 GB of RAM. The models were solved using Gurobi 11.0.0 with default settings, and  a time limit of 3600 seconds.

\subsection{The tolerance model }

We solved the inverse mixed integer optimization  model 
[\text{\bf IMIO($\hat{x},\mathring{c},\tau$)}]
with the tolerance parameter $\tau$ depending on  $\mathring{c}^\top \mathring{x}$  where $\mathring{x}$ is the solution obtained by solving  the forward model with  $\mathring{c}$ for 30 seconds.  It is set to 
    $e${-3} if $\mathring{c}^\top \mathring{x} < e^{3}$,
    $e${-4} if $e^{3} \leq \mathring{c}^\top \mathring{x} < e^{4} $,
    $e$-5  if $e^{4} \leq \mathring{c}^\top \mathring{x} < e^{5} $, and 
    $e${-6}  otherwise.
The resulting cost vector is denoted by $\hat{c}$. To test if it makes $\hat{x}$ optimal, we solve the  forward MIP model for 30 minutes to obtain  an incumbent upper bound ($\overline{ub}$) and a lower bound    
($\overline{lb}$) to compare against.
 To assess the optimality gap of $\hat{x}$ relative  to $\overline{lb}$ and the proximity of $\hat{c}$ to $\mathring{c}$, we look at  their  relative gaps: 
\[\text{RGap} = \frac{|\hat{c}^\top\hat{x}-\overline{lb}|}{\max(1,|\hat{c}^\top\hat{x}|)}
\]
\[
\quad \quad \text{RNorm} = \frac{||\hat{c}||_{1} - ||\mathring{c}||_{1} }{\max(1,||\mathring{c}||_{1})}.\] 
Table \ref{groups_summary2} displays the group number and the number of instances it contains,  
the number of instances in which $\hat{x}$ is optimal under  relative gaps of ($e$-5) and ($e$-2), the minimum (min), maximum (max), and average (avg) relative gaps and  norms, and the computational time in seconds (CPU(s)) for each group. The number of instances in each group is displayed in parentheses and 
\begin{table}[H]
\tabcolsep 5 pt
\centering \small 
\begin{tabular}{r|cc|rrr|rrr|rrr}
\hline
& 
\multicolumn{2}{c|}{\#Optimal} &
\multicolumn{3}{c|}{RGap}&\multicolumn{3}{c|}{RNorm}&\multicolumn{3}{c}{CPU  (s)}\\
\hline
Group (\#instances)&($e$-5)&($e$-2)&min&avg&max&min&avg&max&min&avg&max\\
\hline
Group 1  (18)& 4 & 15 & 0.00 & 0.01 & 0.04 & 0.01 & 0.45 & 1.00 & 0.00 & 0.09 & 0.32 \\
Group 2  (15)& 7 & 12 & 0.00 & 0.02 & 0.14 & 0.08 & 0.68 & 1.00 & 0.01 & 0.04 & 0.13 \\
Group 3  (18)& 9 & 15 & 0.00 & 0.01 & 0.06 & 0.00 & 0.33 & 0.73 & 0.01 & 0.07 & 0.26 \\
Group 4  (18)& 5 & 17 & 0.00 & 0.01 & 0.18 & 0.01 & 0.42 & 0.70 & 0.02 & 0.10 & 0.35 \\
Group 5  (18)& 10 & 18 & 0.00 & 0.00 & 0.00 & 0.01 & 0.37 & 0.68 & 0.02 & 0.11 & 0.22 \\
Group 6  (18)& 6 & 14 & 0.00 & 0.07 & 1.00 & 0.00 & 0.20 & 0.42 & 0.05 & 0.31 & 1.13 \\
Group 7  (18)& 12 & 17 & 0.00 & 0.00 & 0.07 & 0.00 & 0.25 & 0.50 & 0.06 & 0.18 & 0.45 \\
Group 8  (18)& 6 & 17 & 0.00 & 0.00 & 0.01 & 0.00 & 0.58 & 1.00 & 0.06 & 1.08 & 3.90 \\
Group 9  (18)& 10 & 18 & 0.00 & 0.00 & 0.00 & 0.00 & 0.49 & 1.00 & 0.13 & 1.47 & 4.70 \\
Group 10 (18) & 7 & 12 & 0.00 & 0.01 & 0.06 & 0.00 & 0.02 & 0.08 & 0.28 & 0.81 & 2.54 \\
\hline
\end{tabular}
\caption{\small Summary results for the tolerance model:  [\text{\bf IMIO($\hat{x},\mathring{c},\tau$)}].}
\label{groups_summary2}
\end{table}
According to Table \ref{groups_summary2},  the tolerance model solves 76 out of 177 (43\%) and 155 out of 177  (88\%) to an optimality gap of $e$-5  and $e$-2, respectively and leads to very small relative norms, relative gaps and computational times.  The relative gaps, relative norms and computational times are on average
1.3\%, 37.9\% and 0.43 seconds, respectively.
To get a better feel for this remarkable  performance, we examine the detailed results for the largest group, 
 Group 10,  in Table \ref{group10t}.

\begin{table}[H]
\tabcolsep 2.5pt
	\centering
\small 
		\begin{tabular}{r|cccccccccc}
		\multicolumn{11}{l}{}\\ 
			\hline
&		&	&		&	&		&	CPU	&	\multicolumn{2}{c}{Relative} & \multicolumn{2}{c}{Optimal}\\ 
Instance	&	$\overline{ub}$	& $\overline{lb}$	& $\hat{c}^\top \hat{x}$ & $e^\top \varepsilon$	& $||\hat{c}-\mathring{c}||$	&	(s)	&Gap&Norm	&($e$-5)	&	($e$-2)	\\
\hline
air05t1 & 25886.08 & 25885.88 & 25886.08 & 0.22 & 21789.17 & 0.71 & 0.00 & 0.01 & 1 & 1 \\
air05t2 & 26173.90 & 26173.90 & 26174.08 & 0.18 & 17673.21 & 0.69 & 0.00 & 0.00 & 1 & 1 \\
air05t3 & 25927.08 & 25926.86 & 25927.08 & 0.22 & 21912.17 & 0.63 & 0.00 & 0.01 & 1 & 1 \\
neos-3083819-nubut1 & 5564139.66 & 5564137.08 & 5564137.42 & 0.35 & 347430.24 & 0.30 & 0.00 & 0.00 & 1 & 1 \\
neos-3083819-nubut2 & 5480306.94 & 5480302.57 & 5480302.87 & 0.30 & 302221.98 & 0.34 & 0.00 & 0.00 & 1 & 1 \\
neos-3083819-nubut3 & 3561865.64 & 3561534.21 & 3562208.20 & 673.99 & 673994.55 & 0.28 & 0.00 & 0.00 & 0 & 1 \\
neos-1582420t1 & 66.23 & 66.23 & 70.59 & 4.79 & 4790.41 & 0.73 & 0.06 & 0.07 & 0 & 0 \\
neos-1582420t2 & 72.72 & 72.72 & 76.50 & 5.49 & 5487.63 & 0.86 & 0.05 & 0.08 & 0 & 0 \\
neos-1582420t3 & 66.25 & 66.25 & 70.41 & 4.74 & 4737.59 & 0.77 & 0.06 & 0.07 & 0 & 0 \\
nursesched-sprint02t1 & 23.47 & 23.47 & 23.78 & 0.61 & 609.61 & 2.01 & 0.01 & 0.01 & 0 & 0 \\
nursesched-sprint02t2 & 27.96 & 27.96 & 28.52 & 0.58 & 579.46 & 2.54 & 0.02 & 0.01 & 0 & 0 \\
nursesched-sprint02t3 & 26.04 & 26.04 & 26.63 & 0.59 & 585.87 & 1.53 & 0.02 & 0.01 & 0 & 0 \\
drayage-100-23t1 & 11714.44 & 11714.44 & 11714.44 & 0.30 & 296166.55 & 0.45 & 0.00 & 0.02 & 1 & 1 \\
drayage-100-23t2 & 11714.14 & 11714.14 & 11714.42 & 0.28 & 277294.62 & 0.50 & 0.00 & 0.01 & 0 & 1 \\
drayage-100-23t3 & 11714.14 & 11714.14 & 11714.48 & 0.34 & 337123.92 & 0.46 & 0.00 & 0.02 & 0 & 1 \\
drayage-25-23t1 & 11714.14 & 11714.14 & 11714.43 & 0.29 & 292013.61 & 0.80 & 0.00 & 0.02 & 0 & 1 \\
drayage-25-23t2 & 11714.14 & 11714.14 & 11714.42 & 0.28 & 277294.62 & 0.53 & 0.00 & 0.01 & 0 & 1 \\
drayage-25-23t3 & 11714.46 & 11714.46 & 11714.46 & 0.32 & 321442.58 & 0.46 & 0.00 & 0.02 & 1 & 1 \\
\hline\end{tabular}	
\caption{\small Detailed results of the tolerance model, [\text{\bf IMIO($\hat{x},\mathring{c},\tau$)}],
on the largest instances (Group 10).}
\label{group10t}
\end{table}

The columns denote the instance name, the upper bound $\overline{ub}$ and lower bound $\overline{lb}$ when solving the forward problem with the cost vector found ($\hat{c}$), the corresponding objective for  $\hat{c}^\top\hat{x}$, the  LP gap $ e^\top\varepsilon$, the $\ell_1$ norm $||\hat{c}-\mathring{c}||_1$, the relative gap and norm, the CPU time in seconds and  an indicator whether $\hat{x}$ is optimal (1) or not (0) for  relative gaps of $e$-5 and $e$-2.

Table \ref{group10t} reveals that the relative gaps are  very small, making $\hat{x}$  either optimal or very close to optimal. The computational times are less than 2.54 seconds  due to solving an LP. All these features favor the use of the approach for practical  large applications with large data sets. It also  suggests that it can handle additional restrictions  on the cost vector, such as being within a  range and/or requiring it to be integer, as would be  expected in some applications.

\subsection{The bi-objective model} 

For  the bi-objective model [\text{\bf IMIO($\hat{x},\mathring{c},w$)}], we experimented  with different weights and found that 
setting $w_{i} = \max\{\hat{x}_{i},2\}, i \in \cal{I}$ provides good results.  Table \ref{groups_summary3} provides  summary statistics similar to Table \ref{groups_summary2}.

\begin{table}[H]
\tabcolsep 5 pt
\centering \small 
\begin{tabular}{r|cc|rrr|rrr|rrr}
\hline
& 
\multicolumn{2}{c|}{\#Optimal} &
\multicolumn{3}{c|}{RGap}&\multicolumn{3}{c|}{RNorm}&\multicolumn{3}{c}{CPU  (s)}\\
\hline
Group (\#instances)&($e$-5)&($e$-2)&min&avg&max&min&avg&max&min&avg&max\\
\hline
Group 1  (18)& 9 & 16 & 0.00 & 0.00 & 0.01 & 0.02 & 0.38 & 1.00 & 0.00 & 0.04 & 0.12 \\
Group 2  (15)& 9 & 15 & 0.00 & 0.00 & 0.00 & 0.08 & 0.63 & 1.00 & 0.01 & 0.03 & 0.08 \\
Group 3  (18)& 13 & 17 & 0.00 & 0.00 & 0.02 & 0.00 & 0.29 & 0.62 & 0.01 & 0.05 & 0.13 \\
Group 4  (18)& 8 & 10 & 0.00 & 0.05 & 0.32 & 0.09 & 0.43 & 0.67 & 0.02 & 0.05 & 0.11 \\
Group 5  (18)& 9 & 13 & 0.00 & 0.02 & 0.09 & 0.00 & 0.29 & 0.72 & 0.02 & 0.07 & 0.12 \\
Group 6  (18)& 6 & 8 & 0.00 & 0.12 & 1.00 & 0.00 & 0.19 & 0.38 & 0.02 & 0.17 & 0.87 \\
Group 7  (18)& 14 & 14 & 0.00 & 0.09 & 0.76 & 0.00 & 0.23 & 0.50 & 0.04 & 0.13 & 0.40 \\
Group 8  (18)& 12 & 18 & 0.00 & 0.00 & 0.00 & 0.00 & 0.52 & 1.00 & 0.05 & 0.76 & 3.14 \\
Group 9  (18)& 9 & 12 & 0.00 & 0.05 & 0.20 & 0.00 & 0.44 & 1.00 & 0.08 & 0.74 & 2.59 \\
Group 10 (18) & 16 & 16 & 0.00 & 0.00 & 0.02 & 0.00 & 0.02 & 0.08 & 0.20 & 0.41 & 1.06 \\
\hline
\end{tabular}
\caption{\small Summary results for the bi-objective model: [\text{\bf IMIO($\hat{x},\mathring{c},w$)}].}
\label{groups_summary3}
\end{table}

The bi-objective model solves 105 out of 177 (59\%) and 139 out of 177  (79\%) to an optimality gap  of $e$-5  and $e$-2, respectively and leads to very small relative norms, relative gaps and computational times.  The relative gaps, relative norms and computational times are on average 
3.3\%, 34.2\% and 0.25 seconds, respectively.
Similar to Table \ref{group10t},  Table \ref{group10w}, provides detailed results for the largest instances.

\begin{table}[H]
\tabcolsep 1pt
	\centering
\small 
				\begin{tabular}{r|cccccccccc}
		\multicolumn{11}{l}{}\\ 
			\hline
&		&	&		&	&		&	CPU	&	\multicolumn{2}{c}{Relative} & \multicolumn{2}{c}{Optimal}\\ 
Instance	&	$\overline{ub}$	& $\overline{lb}$	& $\hat{c}^\top \hat{x}$ & $e^\top\varepsilon$	& $||\hat{c}-\mathring{c}||$	&	(s)	&Gap&Norm	&($e$-5)	&	($e$-2)	\\
\hline
air05t1 & 25886.08 & 25886.08 & 25886.08 & 0.00 & 21789.39 & 0.32 & 0.00 & 0.01 & 1 & 1 \\
air05t2 & 26133.08 & 26133.08 & 26133.08 & 0.00 & 17673.39 & 0.30 & 0.00 & 0.00 & 1 & 1 \\
air05t3 & 25877.61 & 25877.61 & 25877.61 & 0.00 & 21912.39 & 0.29 & 0.00 & 0.01 & 1 & 1 \\
neos-3083819-nubut1 & 5616510.00 & 5616510.00 & 5616510.00 & 0.00 & 347446.00 & 0.21 & 0.00 & 0.00 & 1 & 1 \\
neos-3083819-nubut2 & 5497276.57 & 5497261.50 & 5497276.57 & 15.07 & 301195.00 & 0.20 & 0.00 & 0.00 & 1 & 1 \\
neos-3083819-nubut3 & 3487309.29 & 3487309.29 & 3487309.29 & 0.00 & 692098.14 & 0.21 & 0.00 & 0.00 & 1 & 1 \\
neos-1582420t1 & 64.56 & 64.56 & 64.56 & 0.23 & 4799.53 & 0.29 & 0.00 & 0.07 & 1 & 1 \\
neos-1582420t2 & 68.08 & 68.08 & 68.08 & 0.88 & 5496.86 & 0.25 & 0.00 & 0.08 & 1 & 1 \\
neos-1582420t3 & 64.66 & 64.66 & 64.66 & 0.00 & 4747.07 & 0.34 & 0.00 & 0.07 & 1 & 1 \\
nursesched-sprint02t1 & 24.17 & 24.17 & 24.67 & 1.00 & 608.83 & 0.94 & 0.02 & 0.01 & 0 & 0 \\
nursesched-sprint02t2 & 26.76 & 26.76 & 27.16 & 0.40 & 579.83 & 1.06 & 0.01 & 0.01 & 0 & 0 \\
nursesched-sprint02t3 & 24.00 & 24.00 & 24.00 & 0.00 & 587.00 & 0.88 & 0.00 & 0.01 & 1 & 1 \\
drayage-100-23t1 & 11714.14 & 11714.14 & 11714.14 & 0.00 & 296166.85 & 0.36 & 0.00 & 0.02 & 1 & 1 \\
drayage-100-23t2 & 11714.14 & 11714.14 & 11714.14 & 0.00 & 277294.89 & 0.35 & 0.00 & 0.01 & 1 & 1 \\
drayage-100-23t3 & 11714.14 & 11714.14 & 11714.14 & 0.00 & 337124.26 & 0.37 & 0.00 & 0.02 & 1 & 1 \\
drayage-25-23t1 & 11714.14 & 11714.14 & 11714.14 & 0.00 & 292013.90 & 0.35 & 0.00 & 0.01 & 1 & 1 \\
drayage-25-23t2 & 11714.14 & 11714.14 & 11714.14 & 0.00 & 277294.89 & 0.34 & 0.00 & 0.01 & 1 & 1 \\
drayage-25-23t3 & 11714.14 & 11714.14 & 11714.14 & 0.00 & 321442.90 & 0.36 & 0.00 & 0.02 & 1 & 1 \\
\hline\end{tabular}	
\caption{\small Detailed  results for the bi-objective model, [\text{\bf IMIO($\hat{x},\mathring{c},w$)}], on  the largest instances  (Group 10).}
\label{group10w}
\end{table}

According to Table \ref{group10w}, the bi-objective model solves all instances to an optimality gap of  $e$-5 except for two instances, and performs better than the  tolerance model. This is not always the case as the tolerance model performs better on other instances, for example Group 9.
It suggests, however, that the two models could be  used concurrently and the better of the two results would be retained. Counting the number of instances where one of the   approaches solves to optimality reveals that 169 out of 177 (95\%) and 120  out of 177 (68\%) are solved to an optimality gap of $e$-5  and $e$-2, respectively.
In the following section, we assess the performance of the two models relative to the literature.

\subsection{Comparison to the literature}

We compare the tolerance and the bi-objective models  
[\text{\bf IMIO($\hat{x},\mathring{c},\tau$)}]
and
[\text{\bf IMIO($\hat{x},\mathring{c},w$)}] to the state-of-the-art trust region cutting plane approach of \cite{BodurChanZhu} who gratefully shared their detailed results,  specifically  the cutting plane with trust region and early stop.

Table \ref{comp_summary} provides a summary of the comparison where the columns are similar to the previous tables. 
\begin{table}[H]
	\tabcolsep 3pt
	\centering \small 
	\begin{tabular}{r|c|c|cc||c|c|cc||c|c|c}
		\hline
		&
		\multicolumn{4}{c|}{[\text{\bf IMIO($\hat{x},\mathring{c},w$)}]}&\multicolumn{4}{c|}{[\text{\bf IMIO($\hat{x},\mathring{c},t$)}]}&\multicolumn{3}{c}{ {\bf CPTR}}\\
		\hline
				&average&average&\multicolumn{2}{c|}{\#Optimal}
			&average&average&\multicolumn{2}{c|}{\#Optimal}
			&average&average& \\
		Group& RNorm & CPU (s) &($e$-5) & ($e$-2) &RNorm & CPU (s)  &($e$-5) & ($e$-2)&RNorm & CPU (s) &  \#Opt \\ \hline
		1 & 0.38 & 0.04 & 9 & 16 & 0.45 & 0.09 & 4 & 15 & 0.34 & 1217.09 & 12 \\
2 & 0.63 & 0.03 & 9 & 15 & 0.68 & 0.04 & 7 & 12 & 0.58 & 797.00 & 12 \\
3 & 0.29 & 0.05 & 13 & 17 & 0.33 & 0.07 & 9 & 15 & 0.12 & 1833.34 & 9 \\
4 & 0.43 & 0.05 & 8 & 10 & 0.42 & 0.10 & 5 & 17 & 0.21 & 2217.74 & 8 \\
5 & 0.29 & 0.07 & 9 & 13 & 0.37 & 0.11 & 10 & 18 & 0.25 & 1617.39 & 12 \\
6 & 0.19 & 0.17 & 6 & 8 & 0.20 & 0.31 & 6 & 14 & 0.16 & 2882.59 & 4 \\
7 & 0.23 & 0.13 & 14 & 14 & 0.25 & 0.18 & 12 & 17 & 0.19 & 1492.91 & 11 \\
8 & 0.52 & 0.76 & 12 & 18 & 0.58 & 1.08 & 6 & 17 & 0.18 & 1851.00 & 12 \\
9 & 0.44 & 0.74 & 9 & 12 & 0.49 & 1.47 & 10 & 18 & 0.34 & 1386.58 & 13 \\
10 & 0.02 & 0.41 & 16 & 16 & 0.02 & 0.81 & 7 & 12 & 0.02 & 2257.18 & 8 \\ \hline
		\multicolumn{3}{r|}{Total:}&105&139&
		\multicolumn{2}{c|}{}&76&155&
		\multicolumn{2}{c|}{}&101 \\
		\hline
	\end{tabular}
	\caption{\small Comparison of  [\text{\bf IMIO($\hat{x},\mathring{c},w$)}], [\text{\bf IMIO($\hat{x},\mathring{c},\tau$)}] and the Trust Region Method of  \cite{BodurChanZhu}}
	\label{comp_summary}
\end{table}

Table \ref{comp_summary} demonstrates that the proposed models provide solution in record computational times as compared to  CPTR. This is due to the inherent advantage of solving LPs as opposed to MIPs. This, however, comes at a cost.  We do not guarantee the optimality of $\hat{x}$ nor the optimality of the norm. Still,  [\text{\bf IMIO($\hat{x},\mathring{c},w$)}]  provides more solutions that render $\hat{x}$ optimal compared to CPTR, but  CPTR results in smaller RNorm values.

One of the advantages of the proposed approach is its  ability to identify cost vectors and solutions that are within the vicinity of  $\hat{x}$. To exploit this property, we perform a limited number of iterations using the classical  cutting plane algorithm. In other words, once a cost vector $\hat{c}$ is identified, we solve $\min \{\hat{c}^\top x: Ax=b, x\geq 0$,  $x_j$ integer $ j=q+1,...,n$  $\}$
for a short time. If $\hat{x}$ is not optimal,   cuts  $\hat{x}^\top c \leq \bar{x}^\top c$ are added for all feasible solutions $\bar{x}$ identified. 

At each iteration, the identified cuts are added and  the tolerance parameter $\tau$  is updated. The cutting plane method  stops when the absolute optimality gap is less than $e$-2, the number of iterations exceeds 1000, or a  time limit of 3600s is reached.
Specifically,  the tolerance model is initially solved with a tolerance of 0.01 to obtain $\hat{c}$, then {\bf [FMIO($c$)]}  is solved for 30 seconds and the resulting cuts are added to [\text{\bf IMIO($\hat{x},\mathring{c},\tau$)}]. The parameter $\tau$ is then re-initialized to 1 and the cutting plane method is started.
At iteration $k$,  $\tau$ is increased to $1.25\tau_{k-1}$ if $k$ is  odd, and to $0.75\tau_{k-1}$  if  $k$ is  even.

We tested on the smallest instances of Group 1, on which CPTR solves the first 12 but struggles with the remaining  6 instances. The results are displayed in Table \ref{CPsummary}. The first column is  
 the final  tolerance $\tau=e^{\kappa}$ attained where only the exponent  is displayed ($\kappa$). Subsequent columns correspond to 
  the objective at $\hat{x}$ ($\hat{c}^\top\hat{x}$),
the objective of the exact MIP solution when solved with  $\hat{c}$ ($\hat{z}$), 
the objective at $\hat{x}$ for the CPTR cost vector $c^{TR}$ ($\hat{x}^\top c^{TR}$),
the objective of the exact MIP solution when solved with  $c^{TR}$ ($z^{TR}$),
the number of  cutting plane iterations performed (itr), 
 the CPU time in seconds and the ratio of ${\|\hat{c}-\mathring{c}\|_1}$ 
to $ {\|c^{TR}-\mathring{c}\|_1}$, respectively.

\begin{table}[H]
	\tabcolsep 1.5pt
	\centering \small 
	\hspace*{-.4cm}
	\begin{tabular}{r|cccccccccc}
		\hline
		Instance	&	$\kappa$	&	$\hat{c}^\top \hat{x}$	&	$\hat{c}^\top x^{IP}$	&	$\hat{x}^\top c^{TR}$	&	$z^{TR}$	&	$\|\hat{c}$-$\mathring{c}\|$	&	$\|c^{TR}-\mathring{c}\|$	&	itr	&	CPU(s)	&	$\frac{\|\hat{c}-\mathring{c}\|}{\|c^{TR}-\mathring{c}\|}$	\\ \hline
markshare-4-0t1 & -0.06 & 0 & 0 & 0 & 0 & 4 & 4 & 5 & 0.12 & 1 \\
markshare-4-0t2 & -0.06 & 0 & 0 & 0 & 0 & 4 & 4 & 5 & 0.1 & 1 \\
markshare-4-0t3 & 0.04 & 0 & 0 & 0 & 0 & 3 & 3 & 6 & 0.42 & 1 \\
gen-ip002t1 & -0.64 & -5592.8 & -5592.8 & -5592.8 & -5592.8 & 198.41 & 198.41 & 47 & 95.64 & 1 \\
gen-ip002t2 & -0.18 & -5885.82 & -5885.82 & -5885.82 & -5885.82 & 170.61 & 170.61 & 22 & 15.85 & 1 \\
gen-ip002t3 & -1.4 & -4644.5 & -4644.5 & -5616.13 & -5616.13 & 617.24 & 249.29 & 101 & 299.57 & 2.48 \\
neos5t1 & -0.87 & 0 & 0 & 0 & 0 & 36 & 36 & 63 & 11.14 & 1 \\
neos5t2 & -0.81 & 0 & 0 & 0 & 0 & 34 & 34 & 59 & 8.45 & 1 \\
neos5t3 & -0.9 & 0 & 0 & 0 & 0 & 36 & 36 & 65 & 14.15 & 1 \\
markshare2t1 & 0.01 & 0 & 0 & 0 & 0 & 7 & 7 & 8 & 0.32 & 1 \\
markshare2t2 & 0.04 & 0 & 0 & 0 & 0 & 7 & 7 & 6 & 0.23 & 1 \\
markshare2t3 & 0.01 & 0 & 0 & 0 & 0 & 7 & 7 & 8 & 0.32 & 1 \\ \hline
pg5\_34t1 & -5.17 & -12237.81 & -12237.82 & -12029.91 & -12063.69 & 6997.9 & 4915.53 & 378 & 2098.18 & 1.42 \\
pg5\_34t2 & -5.82 & -12004.85 & -12004.86 & -12333.3 & -12457.94 & 8081.77 & 5110.66 & 424 & 2331.65 & 1.58 \\
pg5\_34t3 & -4.56 & -12119.29 & -12119.3 & -12170.44 & -12231.1 & 7452.41 & 5200.77 & 334 & 1819.5 & 1.43 \\
pgt1 & -4.42 & -4061.46 & -4061.45 & -5456.09 & -6701.25 & 11100.81 & 5401.63 & 324 & 1226.3 & 2.06 \\
pgt2 & -4.86 & -4136.85 & -4136.66 & -7284.97 & -9287.06 & 11653.14 & 5324.24 & 356 & 1349.03 & 2.19 \\
pgt3 & -4.47 & -4249.19 & -4249.13 & -5208.18 & -6838.91 & 10315.83 & 4666.2 & 328 & 1244.67 & 2.21 \\
		\hline
	\end{tabular}
	\caption{\small
Comparison of  [\text{\bf IMIO($\hat{x},\mathring{c},w$)}]-based cutting plane method and  CPTR on Group 1}
	\label{CPsummary}
\end{table}
Clearly, the approach is capable of finding optimal solutions in competitive computational times. Although the norm $\ell_1$ deviation is not guaranteed, the results show that 
the approach is capable of finding the  optimal solution in all but one instance (gen-ip002t3)
for the instances that the CPTR solved to optimality. For the ones where CPTR can not find an optimal solution within an hour, the approach takes longer time.
In terms of norm, we note that the CPTR norm is a lower bound on the optimal norm, so the norm deviations for the last 6 instances  could very well be optimal as the deviation may be due to the lower bound.

\section{Conclusion and Future Research} \label{sec:con}

 We proposed an inverse optimization framework for mixed integer optimization based on analytic center concepts from interior point methods. The novelty of the idea is centered around the expectation that mixed integer solutions tend to be interior points. The final  model is a linear program that is computationally very efficient. We proposed two variants: the tolerance model and the bi-objective model,  which  achieve remarkable results despite not guaranteeing that the given solution $\hat{x}$ is optimal. In order to guarantee optimality, we complement the tolerance model with a classical cutting plane method based on solving the forward mixed integer program for a short time. This leads to  cost vectors that
are often optimal with respect to the norm, in addition to guaranteeing the optimality of  the given mixed integer solution.  

The proposed approach could benefit from parameter tuning, specifically the tolerance parameter $\tau$ and the objective weight $w_i$. Furthermore, 
the  approach's computational efficiency opens up doors for a multitude of applications where the  objective is inferred from data. These could include contextual objectives and multiple objectives. It would also enable methodological extensions to better handle the trade-off between the two objectives of the bi-level model, the exploration of other nonlinear norm deviations other than the simple $
\ell_1$ norm, and possibly the integration of the approach within search or branching methods to improve the cost deviation from the nominal cost.

\vspace*{.5cm}
\noindent {\bf Acknowledgment:}   We thank  Ian Yihang Zhu, Tim Chan, and Merve Bodur for making the trust region data and results available for testing and comparison. The research was  supported by an NSERC Discovery grant RGPIN-2022-03530 (Elhedhli).







\end{document}